\documentclass[12pt]{article}
\usepackage{amsmath}
\usepackage[pdftex]{graphicx}
\newtheorem{theo}{Theorem}
\newtheorem{lemme}{Lemma}
\newtheorem{question}{Question}
\newtheorem{prop}[lemme]{Proposition}
\newtheorem{cor}[lemme]{Corollary}
\newtheorem{defi}[lemme]{Definition}
\newtheorem{rem}[lemme]{Remark}
\newtheorem{exemple}[lemme]{Example}
\def\R{{\bf R}}

\def\Z{{\bf Z}}

\def\inv#1{\frac{1}{#1}}
\def\qed{\vbox{\hrule height 5pt width 5pt}}
\def\proof{\par\medskip{\bf Proof.}\hskip5pt}
\def\dh{d^{H}}
\def\n#1{|\hskip-1pt |#1|\hskip-1pt |}
\def\fleche#1{\smash{\mathop{\longrightarrow}\limits^{#1}}}

\title{Submanifolds and differential forms on Carnot manifolds, after M. Gromov and M. Rumin}
\author{Pierre Pansu}
\begin{document}
\maketitle

\section{Introduction}

The purpose of these notes is to explain parts of Gromov's survey \cite{Gromov-CC}, in the light of subsequent results of M. Rumin, \cite{Rumin-99}. Among the rich material of \cite{Gromov-CC}, most of which pertains to analysis on metric spaces, we chose to concentrate on the H\"older equivalence problem for Carnot manifolds.

\subsection{Carnot manifolds}

\begin{defi}
\label{defcm}
We shall call \emph{Carnot manifold} the data of a smooth manifold $M$ and a smooth subbundle $H$ of the tangent bundle $TM$, satisfying the \emph{bracket generating condition} : for each $x\in M$, the values at $x$ of iterated Lie brackets of vectorfields tangent to $H$ generate $T_x M$.
\end{defi}
Given a smooth euclidean structure on $H$, the \emph{Carnot-Caratheodory metric} is obtained by minimizing the length of \emph{horizontal curves}, i.e. curves tangent to $H$. The bracket generating condition implies that this distance is finite.

\medskip

\textbf{Question.} How far can a Carnot-Caratheodory metric be from a Riemannian metric ?

\begin{exemple}
\label{exheis}
Consider the group $Heis^3$ of real unipotent $3\times 3$ matrices $\begin{pmatrix}
1&x&z \\ 0&1&y\\ 0&0&1
\end{pmatrix}$. The differential 1-form $\theta=dz-xdy$ on $Heis^3$ is left invariant. So is its kernel $H=ker(\theta)$.
\end{exemple}
One can use the left invariant Euclidean structure $dx^2 +dy^2$ to define a Carnot-Caratheodory metric. At small scales, this metric is very different from any Riemannian metric in 3 dimensions. Indeed, its Hausdorff dimension is 4 instead of 3. This follows from the existence of the one parameter group of homothetic diffeomorphisms
\begin{eqnarray*}
\delta_\epsilon : (x,y,z)\mapsto (\epsilon x,\epsilon y ,\epsilon^2 z).
\end{eqnarray*}
Since $\delta_\epsilon$ takes unit balls to balls of radius $\epsilon$ and multiplies volumes by $\epsilon^4$, one needs $\epsilon^{-4}$ $\epsilon$-balls to cover a bounded  open set.

\subsection{Carnot groups}
The above Hausdorff dimension calculation immediately extends to the following family of examples.

\begin{defi}
\label{defcg}
A Carnot group is a simply connected Lie group $G$ equipped with a subspace $V^1 \subset Lie(G)$ such that the subspaces defined inductively by $V^i =[V^1 ,V^{i-1}]$ constitute a \emph{gradation} of $Lie(G)$, i.e.
\begin{eqnarray*}
Lie(G)=V^1 \oplus \cdots \oplus V^r , \quad \textrm{and}\quad [V^{i},V^{j}]\subset V^{i+j}.
\end{eqnarray*}
\end{defi}
Left translating $V^1$ yields a subbundle $H$ which satisfies the bracket generating condition. Choose a Euclidean structure on $V^1$ and left-translate it. The group automorphisms defined on the Lie algebra by
\begin{eqnarray*}
\delta_\epsilon (v)=\epsilon^i v \quad\textrm{for}\quad v\in V^i ,
\end{eqnarray*}
are homothetic for the Carnot-Caratheodory metric. It follows that the Hausdorff dimension of this metric is equal to
\begin{eqnarray*}
\sum_{i=1}^{r} i\,dim(V^i ).
\end{eqnarray*}

\subsection{Tangent cones}
Carnot groups play, in the family of Carnot manifolds, the role played by Euclidean space among Riemannian manifolds, at least under some restrictive condition.

\begin{defi}
\label{defer}
Let $(M,H)$ be a Carnot manifold. For $x\in M$, define $H^2 (x)$ as the linear span of values at $x$ of brackets of vector fields tangent to $H$. And recursively, let $H^i (x)$ be the linear span of values at $x$ of brackets of sections of $H$ and of $H^{i-1}$. Say $H$ is \emph{equiregular} if $x\mapsto dim(H^i (x))$ is constant for all $i$.
\end{defi}

\begin{exemple}
\label{exner}
The kernel $H$ of the differential 1-form $dz-x^2 dy$ on $\R^3$, known as the \emph{Martinet} Carnot structure, is not equiregular.
\end{exemple}
Indeed, when $x\not=0$, the generating vectorfields $\partial_x$ and $\partial_y +x^2 \partial_z$ and their Lie bracket $2x\partial_z$ are linearly independent, so that $H^2 =\R^3$. At points where $x=0$, $H^2 =H$. Still, $H^3=\R^3$, thus the bracket generating condition is satisfied.

\begin{exemple}
\label{exer}
In 3 dimensions, equiregular Carnot manifolds coincide with contact manifolds. 
\end{exemple}

\begin{theo}
\label{NSWM}
\emph{(Nagel-Stein-Wainger \cite{NSW}, Mitchell \cite{Mitchell}).} An equiregular Carnot manifold is asymptotic, at each point $x$, to a Carnot group $G_x$ called its \emph{tangent cone} at $x$. It follows that
\begin{eqnarray*}
\textrm{Hausdorff dimension}=Q=:\sum_{i=1}^r i\,(dim(H^i )-dim(H^{i-1})).
\end{eqnarray*}
\end{theo}

It follows for instance that equiregular Carnot manifolds are never biLipschitz homeomorphic to Riemannian manifolds.

\subsection{BiLipschitz equivalence}

\begin{theo}
\label{thmp}
\emph{(P. Pansu, \cite{Pansu-89}, see also \cite{Vodopyanov}).} Two Carnot groups are biLipschitz homeomorphic (resp. quasiconformally homeomorphic) if and only if they are isomorphic.
\end{theo}

\begin{theo}
\label{thmmm}
\emph{(G. Margulis, G. Mostow, \cite{MM1}, \cite{MM2}, see also \cite{IV}).} If $f:M\to M'$ is a quasiconformal homeomorphism of equiregular Carnot manifolds, then for all $x\in M$, $G'_{f(x)}$ is isomorphic to $G_x$.
\end{theo}

\begin{question}
\label{qlip}
Assume two equiregular Carnot manifolds $M$ and $M'$ are quasiconformally homeomorphic. Does there exist a diffeomorphism $M\to M'$ mapping $H$ to $H'$ ?
\end{question}

\begin{question}
\label{qliph}
Assume an equiregular Carnot manifold $M$ is quasiconformally homogeneous, i.e. for every pair of points $x$, $x'\in M$, there exists a quasiconformal homeomorphism of $M$ mapping $x$ to $x'$. Does there exist a transitive $H$-preserving action of a finite dimensional Lie group on $M$ ?
\end{question}

\subsection{H\"older equivalence}

Since the biLipschitz equivalence problem seems to be understood to some extent, we turn to a harder problem : when are Carnot manifolds H\"older equivalent ?

\begin{theo}
\label{thmrc}
\emph{(Rashevski \cite{Rashevski}, Chow \cite{Chow},...).} Let $(M,H)$ be a Carnot manifold with $H^r =TM$. Let $g$ be a Riemannian metric on $M$. Then identity $(M,g)\to(M,H)$ is locally of class $C^{1/r}$ and its inverse is locally Lipschitz.
\end{theo}

\begin{rem}
\label{remholder}
Let $(M,h)$ be an equiregular Carnot manifold of dimension $n$ and Hausdorff dimension $Q$. If $\alpha>n/Q$, there are no $\alpha$-H\"older-continuous homeomorphisms of Riemannian manifolds to $M$.
\end{rem}
Indeed, if $M'$ is $n$-dimensional Riemannian, then $dim_{Hau}M'=n$. If $f(M')\subset M$ is open, $dim_{Hau}f(M')=Q$, thus the following lemma implies that $\alpha\leq n/Q$. 

\begin{lemme}
\label{dimholder}
If $f:M' \to M$ is a homeomorphism which is $C^{\alpha}$-H\"older continuous, then
\begin{eqnarray*}
\alpha \,dim_{Hau}f(M')\leq dim_{Hau}M'.
\end{eqnarray*}
\end{lemme}

\begin{defi}
\label{defalpha}
Let $\alpha(M,H)$ be the supremum of exponents $\alpha$ such that there (locally) exists a $\alpha$-H\"older-continuous homeomorphism of $\R^n$ onto an open subset of $M$.
\end{defi}

\begin{exemple}
\label{alphaheis}
Theorem \ref{thmrc} and remark \ref{remholder} imply that $1/2 \leq \alpha(Heis^3 )\leq 3/4$.
\end{exemple}

\begin{question}
\label{bestalpha}
Find estimates for $\alpha(M,H)$. For instance, is it true that $\alpha(Heis^3 )=1/2$ ?
\end{question}

Not much is known. For instance, the best known upper bound for $\alpha(Heis^3 )$ is $2/3$, which, as we shall see, follows from the isoperimetric inequality.

\subsection{Results to be covered}

Following Gromov, \cite{Gromov-CC}, we shall give two proofs of the isoperimetric inequality in Carnot manifolds. 

The first one relies on the wealth of horizontal curves. More generally, again following Gromov, \cite{Gromov-PDR}, we shall show that certain Carnot manifolds admit plenty of horizontal $k$-dimensional manifolds, which can be used to prove that $n-k$-dimensional topological manifolds have Hausdorff dimension $\geq Q-k$. This allows to sharpen the upper bound on $\alpha(M,H)$ given by the isoperimetric inequality.

The second one makes a clever use of differential forms. We shall describe a deformation of the de Rham complex of a Carnot manifold, discovered by M. Rumin, \cite{Rumin-99}, which gives alternative proofs of upper bounds on $\alpha(M,H)$. In fact, a combination of Gromov's and Rumin's ideas provides bounds in terms of the homology of the tangent cone which are rather easily computable for every given Carnot group, see Corollary \ref{corweight}, and cover all known results. Unfortunately, these bounds are never sharp.

\subsection{Acknowledgements}

Thanks to Ya. Eliashberg and M. Rumin for sharing their understanding of the subject, and to D. Isangulova for carefully reading the manuscript.
 
\section{Hausdorff dimension of hypersurfaces}

\subsection{The isoperimetric inequality}

Let $(M,H)$ be an equiregular Carnot manifold of Hausdorff dimension $Q$. For simplicity, we denote by $vol$ the $Q$-dimensional Hausdorff measure, and by $area$ the $Q-1$-dimensional Hausdorff measure. The following inequality is due to N. Varopoulos, \cite{Varopoulos}, in the case of Carnot groups (see also \cite{Pansu-82} for the case of $Heis^3$), with a rather sophisticated proof. 

\begin{theo}
\label{thmi}
Let $K$ be a compact subset in an equiregular Carnot manifold of Hausdorff dimension $Q$. There exist constants $c$ and $C$ such that for every piecewise smooth domain $D\subset K$, 
\begin{eqnarray*}
vol(D)\leq c \Rightarrow vol(D)\leq C\,area(\partial D)^{Q/Q-1}.
\end{eqnarray*} 
\end{theo}

\begin{rem}
\label{inenonreg}
Gromov, \cite{Gromov-CC} page 166, observes that the proof applies as well to non equiregular Carnot manifolds, provided the definition of $area$ be adapted. 
\end{rem}

\begin{rem}
\label{ine}
In case $(M,H)$ is a Carnot group, the inequality is valid for arbitrary relatively compact open sets.
\end{rem}
Indeed, by dilation homogeneity, the constants do not depend on the compact set $K$.

\begin{cor}
\label{idimalpha}
Let $(M,H)$ be an equiregular Carnot manifold of dimension $n$ and Hausdorff dimension $Q$. Then
\begin{eqnarray*}
\alpha(M,H)\leq \frac{n-1}{Q-1}.
\end{eqnarray*}
\end{cor}

\subsection{Flow tube estimate}

Gromov's proof (\cite{Gromov-CC}, pages 159-164) relies on pretty general principles. 

Given a vectorfield $X$ on $M$ with (locally defined) flow $\phi_t$ and a subset $S\subset M$, let $Tube(S,\tau)$, the \emph{tube} on $S$ be
\begin{eqnarray*}
Tube(S,\tau)=\{\phi_t (s)\,|\,s\in S,\,0\leq 0\leq \tau\}.
\end{eqnarray*}

\begin{lemme}
\label{tube}
Let $X$ be a horizontal vectorfield on $M$. Let $\epsilon$ and $\tau$ be small (depending on $K$). Let $B$ be an $\epsilon$-ball such that $Tube(B,\tau)$ is contained in $K$. Then
\begin{eqnarray*}
vol(Tube(B,\tau))\leq \textrm{const.}\,\frac{\tau}{\epsilon}vol(B). 
\end{eqnarray*}
It follows that for arbitrary $S\subset K$ such that $Tube(S,\tau)\subset K$,
\begin{eqnarray*}
vol(Tube(S,\tau))\leq \textrm{const.}\,\tau \,area(S),
\end{eqnarray*}
where the constant depends only on $X$ and on $K$.
\end{lemme}

\proof
According to \cite{NSW}, one can choose coordinates such that $X=\frac{\partial}{\partial x_1}$ and $B$ is contained in a box $\{\forall i,\,|x_i |\leq \epsilon^{w(i)}\}$ of volume $\sim\epsilon^{\sum w(i)}\leq\textrm{const.}vol(B)$. Then 
\begin{eqnarray*}
Tube(B,\tau)\subset \{-\epsilon\leq x_1 \leq\tau+\epsilon\,\textrm{and}\,\forall i\geq 2,\,|x_i |\leq \epsilon^{w(i)}\}, 
\end{eqnarray*}
thus 
\begin{eqnarray*}
vol(Tube(B,\tau))\leq \textrm{const.}\tau \,\epsilon^{\sum_{i\geq 2} w(i)}\leq \textrm{const.}\frac{\tau}{\epsilon}vol(B). 
\end{eqnarray*}
Cover set $S$ with small balls $B_j$ with radii $r_j$. There exists a constant $\eta(K)$ such that $\eta\,vol(B_j )\leq r_{j}^{Q}$. Then $Tube(S,\tau)\subset \bigcup_j Tube(B_j ,\tau)$, thus
\begin{eqnarray*}
\sum_{j}r_{j}^{Q-1}
&\geq&
\eta\sum_{j}r_{j}^{-1}vol(B_j )\\
&\geq&
\textrm{const.}\tau^{-1}\sum_{j}vol(Tube(B_j ,\tau))\\
&\geq&
\textrm{const.}\tau^{-1} vol(Tube(S,\tau)).\qed
\end{eqnarray*}

\subsection{Local isoperimetric inequality}

Given smooth vectorfields $X_1 ,\ldots,X_k$, denote by $Tube_j (S,\tau)$ the $\tau$-tube generated by $X_j$ and by $MTube(S,\tau)$ the \emph{multitube} i.e. the set of points reached, starting from a point in $S$, by flowing $X_1$ during some time $t_1 \leq \tau$, then flowing $X_2$ for time $t_2 \leq \tau$, and so on. 

The bracket generating condition allows to choose smooth horizontal vectorfields $X_1 ,\ldots,X_k$ such that the $\tau$-multitube of any point $x\in K$ under them contains $B(x,\tau)$, for all $\tau\leq \textrm{const.}(K)$. Then there exists a constant $\lambda=\lambda(K)\geq 1$ such that for all $x\in K$,
\begin{eqnarray*}
MTube(\{x\},\tau)\subset B(x,\lambda\tau).
\end{eqnarray*}  

\begin{prop}
\label{lii}
For every ball $B$ of radius $R\leq\textrm{const.}(K)$, such that the concentric ball $\lambda B\subset K$, and for every open subset $D\subset K$ with $vol(D\cap \lambda B)\leq \inv2 vol(B)$,
\begin{eqnarray*}
vol(D\cap B)\leq \textrm{const.}(K)\,R\,area((\partial D) \cap \lambda B).
\end{eqnarray*}
\end{prop}

\proof
Let $\tau=2R$. Some significant portion $D_0$ of $D\cap B$ must be carried out of $D$ by the flow of some vectorfield $X_i$. Indeed, otherwise, the multitube $MTube(D\cap B,\tau)$ would be almost entirely contained in $D$. But this multitube contains $B$ which has volume at least twice that of $D$, contradiction\footnote{Here, I merely copy without understanding Gromov's one sentence proof the trivial Measure Moving Lemma}{}. Since the multitube is entirely contained in $\lambda B$, Lemma \ref{tube} applied to $X_i$ then gives 
$$
vol(D_0 )\leq vol(T_i (D\cap B,\tau))\leq const.\,R\,area(\partial D\cap\lambda B).\qed
$$

\subsection{Covering Lemma}

\begin{lemme}
\label{cov}
If $vol(D)\leq\textrm{const.}(K)$, there exists a collection of disjoint balls $B_{j}$ such that
\begin{itemize}
\item $D$ is covered by concentric balls $2B_{j}$.
\item $vol(D\cap\lambda^{-1}B_{j})\geq \frac{1}{2}vol(\lambda^{-2}B_{j})$.
\item $vol(D\cap B_{j})\leq \frac{1}{2}vol(\lambda^{-1}B_{j})$.
\end{itemize}
\end{lemme}

\proof
Fix a radius $R=R(K)$ such that all $R$-balls contained in $K$ have roughly the same volume $\textrm{const.}R^Q$. Fix $\textrm{const.}(K)$ such that if $vol(D)\leq\textrm{const.}(K)$, then
\begin{eqnarray*}
\frac{vol(D)}{\textrm{const.}(\lambda^{-1}R)^Q}\leq \frac{1}{2}.
\end{eqnarray*}
Then for all $x\in D$,
\begin{eqnarray*}
\frac{vol(D\cap B(x,R))}{vol(\lambda^{-1}B(x,R))}\leq \frac{1}{2}.
\end{eqnarray*}

Given $x\in D$, consider the sequence of concentric balls $\beta_{\ell}=B(x,\lambda^{-\ell}R)$. Since $D$ is open, the ratio 
\begin{eqnarray*}
\frac{vol(D\cap \beta_{\ell})}{vol(\lambda^{-1}\beta_{\ell})}
\end{eqnarray*}
tends to $\lambda^{Q}\geq 1$. Let $B(x)$ be the last ball in the sequence $\beta_{\ell}$ such that this ratio is less than $\frac{1}{2}$. By construction, the balls $B(x)_{x\in D}$ cover $D$ and satisfy two of the assumptions of the lemma.

Order the balls $B(x)$ according to their radii, pick a largest one, call it $B_{0}$, then pick a largest one among those which do not intersect $B_{0}$, call it $B_{1}$, and so on. In this way, one obtains a collection of disjoint balls. If $x\in D$ and $B(x)$ has not been selected, then $B(x)$ intersects some selected ball $B_{j}$ which is larger than $B(x)$. This implies that $x\in 2B_{j}$. Therefore the concentric balls $2B_{j}$ cover $D$.\qed

\subsection{From local to global}

Let $D\subset K$ have $vol(D)\leq\textrm{const.}(K)$. Apply Lemma \ref{cov} to get disjoint balls $B_j$ such that $vol(D\cap B_j )\leq \frac{1}{2}vol(\lambda^{-1}B_j )$. The local isoperimetric inequality \ref{lii} applies in each $B_j$ and yields
\begin{eqnarray*}
vol(D\cap \lambda^{-1}B_j )\leq \textrm{const.}\, R_j \,area((\partial D)\cap B_j ).
\end{eqnarray*} 
Since 
\begin{eqnarray*}
R_{j}^{Q}\leq \textrm{const.}\, vol(\lambda^{-2}B_j )\leq \textrm{const.}\, vol(D\cap\lambda^{-1}B_j ),
\end{eqnarray*}
\begin{eqnarray*}
vol(D\cap \lambda^{-1}B_j )\leq \textrm{const.}\, area((\partial D)\cap B_j )^{Q/Q-1}.
\end{eqnarray*}
Finally,
\begin{eqnarray*}
vol(D\cap 2B_j )&\leq& vol(2B_j )\\
&\leq&\textrm{const.}\, vol(\lambda^{-2}B_j )\\
&\leq&\textrm{const.}\, vol(D\cap \lambda^{-1}B_j )\\
&\leq& \textrm{const.}\, area((\partial D)\cap B_j )^{Q/Q-1}.
\end{eqnarray*}
Since $2B_j$ cover and $B_j$ are disjoint, one can sum up,
\begin{eqnarray*}
vol(D)&\leq& \sum_{j}vol(D\cap 2B_j )\\
&\leq&\textrm{const.}\, \sum_{j}area((\partial D)\cap B_j )^{Q/Q-1}\\
&\leq&\textrm{const.}\, (\sum_{j}area((\partial D)\cap B_j ))^{Q/Q-1}\\
&\leq& \textrm{const.}\, area(\partial D)^{Q/Q-1},
\end{eqnarray*}
where a convexity inequality has been used.

This completes the proof of Theorem \ref{thmi}.\qed

\subsection{Link with Sobolev and Poincar\'e inequalities}

It is a classical fact that isoperimetric inequalities are equivalent to Sobolev type inequalities. For a smooth function $u$ on $M$, let $\dh u$ denote the restriction to $H$ of the differential $du$. 

\begin{prop}
\label{sobolev}
The isoperimetric inequality \ref{thmi} is equivalent to the following \emph{Sobolev inequality}, with the same constant. For all smooth functions $u$ with support contained in $K$,
\begin{eqnarray*}
\n{u}_{Q/Q-1}\leq\textrm{const.}\n{\dh u}_1 .
\end{eqnarray*} 
\end{prop}

\proof
It relies on the coarea formula : for smooth $u:K\to \R$ and positive $f$,
\begin{eqnarray*}
\int_{M}f\,d\mathcal{H}^Q=\int_{\R}(\int_{\{u=t\}}\frac{f}{|\dh u|}\,d\mathcal{H}^{Q-1})\,dt.
\end{eqnarray*}

Assume isoperimetric inequality. Write
\begin{eqnarray*}
|u|=\int_{0}^{+\infty}1_{\{|u|>t\}}\,dt
\end{eqnarray*}
as a sum of characteristic functions. Then take $L^{Q/Q-1}$ norms,
\begin{eqnarray*}
\n{u}_{Q/Q-1}&\leq &\int_{0}^{+\infty}\n{1_{\{|u|>t\}}}_{Q/Q-1}\,dt\\
&=&\int_{0}^{+\infty}vol(\{|u|>t\})^{Q-1/Q}\,dt\\
&\leq&\textrm{const.}\,\int_{0}^{+\infty}area(\{|u|=t\})\,dt\\
&=&\n{\dh u}_1 .
 \end{eqnarray*}

Conversely, apply Sobolev inequality to steep functions of the distance to an open set $D$. This gives back the isoperimetric inequality.\qed

\begin{prop}
\label{poincare}
A slightly strengthened form of the local isoperimetric inequality \ref{lii}, namely 
\begin{eqnarray*}
vol(D\cap B)\leq \inv2 vol(B) \Rightarrow
vol(D\cap B)\leq \textrm{const.}(K)\,R\,area((\partial D) \cap \lambda B),
\end{eqnarray*}
is equivalent to the following \emph{$(1,1)$-Poincar\'e inequality}. For a smooth function $u$ defined on a ball $\lambda B$ of radius $R$,
\begin{eqnarray*}
\inf_{c\in\R}\int_{B}|u-c|
\leq \textrm{const.}\,R\,\int_{\lambda B}|\dh u|.
\end{eqnarray*}
\end{prop}

\proof
Assume local isoperimetric inequality holds. Up to replacing $u$ with $u-c$ for some constant $c$, one can assume that
\begin{eqnarray*}
vol(\{u>0\}\cap B)\leq \frac{1}{2}vol(B),\quad vol(\{u<0\}\cap B)\leq \frac{1}{2}vol(B).
\end{eqnarray*}
Write $u=u_+ -u_-$ where $u_+ =\max\{u,0\}$. Then
\begin{eqnarray*}
\int_{B}u_+ &=&\int_{B}(\int_{0}^{+\infty}1_{\{u>t\}}\,dt)\\
&=&\int_{0}^{+\infty}vol(\{u>t\}\cap B)\,dt\\
&\leq&\textrm{const.}\,R\,\int_{0}^{+\infty}area(\{u=t\}\cap\lambda B)\,dt\\
&=&\textrm{const.}\,R\,\int_{\lambda B}|\dh u_+|,
\end{eqnarray*}
by coarea formula. Then
\begin{eqnarray*}
\int_{B}|u|\leq\textrm{const.}\,R\,(\int_{\lambda B}|\dh u_+|+\int_{\lambda B}|\dh u_-|)=\textrm{const.}\,R\,\int_{\lambda B}|\dh u|.
\end{eqnarray*}

Conversely, apply $(1,1)$-Poincar\'e inequality to steep functions of the distance to an open set $D\subset\lambda B$. This gives back the local isoperimetric inequality.\qed

\section{Hausdorff dimension of higher codimensional submanifolds}

Let $M$ be a Carnot manifold. According to lemma \ref{dimholder}, if we can show that all subsets $V\subset M$ of topological dimension $k$ have Hausdorff dimension at least $d$, then $\alpha(M)\leq k/d$. In this section, we prove results of this kind, which in some cases improve on the upper bound obtained in the previous section.

\subsection{A topological criterion}

\begin{prop}
\label{top}
Let $M$ be an $n$-dimensional manifold and $V\subset M$ a subset of topological dimension $dim_{top}(V)\geq n-k$. Then there exists a $k$-dimensional polyhedron $P$ and a continuous map $f:P\to M$ such that any map $\tilde{f}:P\to M$ sufficiently $C^0$-close to $f$ hits $V$. Such a map is called a \emph{transversal} to $V$.
\end{prop}

\proof
This follows from a homological criterion due to Alexandrov, see \cite{Nagata} page 248.\qed

\subsection{Wealth}

We are looking for horizontal immersions in Carnot manifolds $(M,H)$, i.e. immersions whose image is tangent to $H$. We want enough of them to foliate open sets. In view of Lemma \ref{top}, we need to approximate continuous maps from arbitrary polyhedra with immersions. Not every polyhedron is homotopy equivalent to a manifold. Therefore, we enlarge the class of manifolds, by considering spaces, called \emph{branched manifolds} obtained by gluing manifolds along open sets. Given continuous maps $f_0 :P\to M$ and $f:W\to M$, we say that $f$ is $\epsilon$-close to $f_0$ if there exist homotopy equivalences $\phi:W\to P$ and $\phi':P\to W$ such that $\sup |f -f_0 \circ \phi|<\epsilon$, $\sup |\phi'\circ\phi -id_W |<\epsilon$ and $\sup |\phi\circ \phi' -id_P |<\epsilon$.

\begin{center}
\includegraphics[width=2in]{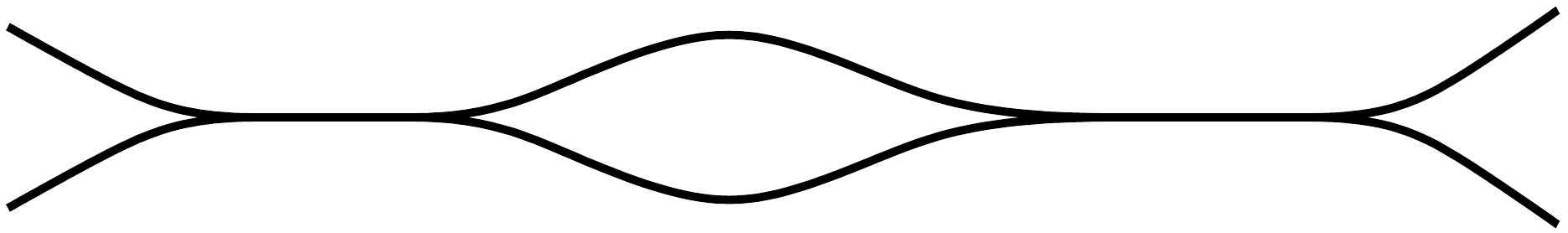}
\par
An immersion of a branched manifold
\end{center}

\begin{defi}
\label{deffoliated}
A \emph{foliated horizontal immersion} in a Carnot manifold $(M,H)$ is a smooth immersion $f:W\times\R^{\ell}\to M$, where $W$ is a branched manifold, and such that for each $z\in\R^{\ell}$, $f_{|W\times\{z\}}$ is horizontal, i.e. tangent to $H$.
\end{defi}

\begin{defi}
\label{defwealth}
Say an $n$-dimensional Carnot manifold $(M,H)$ is \emph{$k$-rich} at a point $m\in M$ if there exists a neighborhood $U$ of $m$ such that, given $\epsilon>0$ and a continuous map from a $k$-dimensional polyhedron $f_0 :P\to U$, there exists a foliated horizontal immersion $f:W\times\R^{n-k}\to M$ with $dim(W)=k$ which is $\epsilon$-close to $f_0$.
\end{defi}

\begin{lemme}
\label{dimrich}
Let $M$ be an $n$-dimensional Carnot manifold. Assume $M$ is $k$-rich at some point $m\in M$. Then for every $n-k$-dimensional subset $V\subset M$ passing through $m$,
\begin{eqnarray*}
dim_{Hau}(V)-dim_{top}(V)\geq dim_{Hau}(M)-dim_{top}(M).
\end{eqnarray*}
If follows that 
\begin{eqnarray*}
\alpha(M)\leq \frac{n-k}{d-k}, \quad d=dim_{Hau}(M).
\end{eqnarray*}\end{lemme}

\proof
Given a foliated horizontal immersion $f:W\times\R^q$ which is close to a transversal to $V$, define the \emph{$\tau$-tube} on $S\subset M$ as the union of $f(B(w,\tau)\times\{z\})$ for which $f(w,z)\in S$. Since $f$ restricted to $W$ factors is horizontal, the $\tau$-tube on an $\epsilon$-ball has volume at most $\textrm{const.}\tau^k \epsilon^{d-k}$. 

Cover $V$ with $\epsilon_j$-balls. The corresponding $\tau$-tubes $T_j$ cover the $\tau$-tube $U$ on $V$. Since $f$ is close to a transversal to $V$, $U$ contains an open set. Then
\begin{eqnarray*}
\sum_{j}\epsilon_{j}^{d-k}\geq \textrm{const.}\tau^{-k}vol(\bigcup_{j}T_j )\geq \textrm{const.}vol(U)
\end{eqnarray*}
is bounded from below, which shows that $dim_{Hau}(V)\geq d-k$.\qed

\subsection{Main result}

\begin{theo}
\label{generich}
\emph{(M. Gromov).} 
A contact structure in dimension $n=2m+1$ is $m$-rich at all points.

Let $0\leq k\leq h\leq n$. Assume that
\begin{eqnarray*}
h-k\geq (n-h)k.
\end{eqnarray*}
Then a generic $h$-dimensional distribution on a $n$-dimensional manifold is $k$-rich at almost every point.
\end{theo}

\begin{cor}
\label{geneholder}
If $(M,H)$ is a $2m+1$-dimensional contact manifold, then $\alpha(M,H)\leq \frac{m+1}{m+2}$.

If $(M,H)$ is a generic Carnot manifold of dimension $n$, Hausdorff dimension $Q$, with $dim(H)=h$ and $h-k\geq (n-h)k$, then $\alpha(M,H)\leq \frac{n-k}{Q-k}$.
\end{cor}

\proof
The proof has three steps.
\begin{enumerate}
  \item Linear algebra : analyze the differential of the equation for horizontal immersions.
  \item Analysis : an implicit function theorem (J. Nash) yields local existence of regular horizontal immersions.
  \item Topology : passing from local to global existence (S. Smale).
\end{enumerate}

In the next three sections, we outline some of the ideas in this proof, following \cite{Gromov-PDR}. An alternate approach to the third step is described in \cite{Gromov-Oka}.

\section{Linearizing horizontality}

\subsection{Isotropic subspaces}

Locally, a $h$-dimensional plane distribution $H$ on a $n$-dimensional manifold $M$ can be viewed as the kernel of a $\R^{n-h}$-valued 1-form $\theta$. 

An immersion $f:V\subset M$ is horizontal iff $f^*\theta=0$. Observe that this implies that $f^*d\theta=0$. 

\begin{defi}
\label{defisotropic}
Let $m\in M$. A linear subspace $S\subset H_m$ is \emph{isotropic} if $d\theta_{|S}=0$. 
\end{defi}

\begin{exemple}
\label{exisotropic}
1-dimensional subspaces are always isotropic. If $H$ is a contact structure on $M^{2m+1}$ (resp. quaternionic contact structure on the sphere $S^{4m+3}$), all isotropic subspaces have dimension $\leq m$.  
\end{exemple}
Here, $S^{4m+3}$ is the unit sphere in the quaternion vector space $\mathbf{H}^{m+1}$, and for $m\in S^{4m+3}$, $H_m$ is the quaternionic hyperplane orthogonal to $m$.
 
\medskip

In particular, a contact (resp. quaternionic contact) manifold has no horizontal immersions of dimension $k>m$.

\subsection{Regular isotropic subspaces}

Our goal is to solve the horizontal immersion equation $E(f)=0$, where, for an immersion $f:V\to M$, $E(f)$ is the $\R^{n-h}$-valued 1-form on $V$ defined by $E(f)=f^*\theta$. 

Let $X$ be a vectorfield along $f$ (i.e. a section of $f^*TM$ on $V$), viewed as a tangent vector at $f:V\to M$ to the space of immersions. The directional derivative of $E$ at $f$ in the direction $X$ is given by Cartan's formula
\begin{eqnarray*}
D_f E(X)=\mathcal{L}_X \theta =d(\iota_X \theta)+f^* (\iota_X (d\theta)).
\end{eqnarray*}

Observe that if $X$ is horizontal, the first term vanishes, and $D_f E(X)$ does not involve any derivatives of $X$. If the second term is pointwise onto, we have an easy way to (right-)invert the operator $D_f E$.

\begin{defi}
\label{defregular}
Let $H=ker(\theta)$ where $\theta$ is $\R^{n-h}$-valued. Let $m\in M$. Say a linear subspace $S\subset H_m$ is \emph{regular} if the linear map
\begin{eqnarray*}
H_m \to Hom(S,\R^{n-h}),\quad X\mapsto (\iota_X d\theta)_{|S}
\end{eqnarray*}
is onto.
\end{defi} 

\begin{exemple}
\label{exregular}
In contact manifolds (resp. on the quaternionic contact $4m+3$-sphere), all isotropic subspaces of $H$ are regular.

A distribution $H$ is sometimes said to satisfy the \emph{strong bracket generating hypothesis} if all 1-dimensional subspaces are regular. Such distributions are very rare.
\end{exemple}
Note that if $H$ admits a regular $k$-dimensional subspace, then $h\geq(n-h)k$. In other words, the codimension of $H$ is very small, $n-h\leq \frac{n}{k+1}$, which is already somewhat restrictive.

\subsection{Genericity of regularity}

The existence of regular isotropic subspaces imposes an even stronger restriction on dimensions. This dimensional condition is genericly sufficient.

\begin{prop}
\label{isoreggeneric}
Let $H$ be a $h$-dimensional distribution on a $n$-dimensional manifold. If $H$ admits a regular $k$-dimensional isotropic subspace at some point, then $h-k\geq(n-h)k$. 

Conversely, if $h-k\geq(n-h)k$, a generic $h$-dimensional distribution admits regular isotropic $k$-planes at almost every point.
\end{prop}

\proof
If $S\subset H_m$ is isotropic, then $S$ is contained in the kernel of the map $X\mapsto \iota_X (d\theta)_{|S}$. If $S$ is regular, this map is onto with a $\geq k$-dimensional kernel. This implies $h-k\geq(n-h)k$. 

Conversely, observe that regular isotropic $k$-planes are the smooth points of the variety of isotropic $k$-planes. Their existence is a Zariski open condition on a 2-form $\omega$. 

The dimension condition guarantees that this open set $U\subset (\Lambda^2 \R^h )\otimes\R^{n-h}$ is non empty. Indeed, pick any surjective linear map $L:\R^{h-k}\to Hom(\R^k ,\R^{n-h})$, viewed as a $h-k\times k$ matrix with entries in $\R^{n-h}$, and let $\omega$ be the 2-form on $\R^k \oplus \R^{h-k}$ with matrix $\begin{pmatrix}
0&L \\ -L^{\top}&0
\end{pmatrix}$. By construction, $\R^k$ is isotropic and the associated map $\R^{h}\to Hom(\R^k ,\R^{n-h})$ is $\begin{pmatrix}
0&L 
\end{pmatrix}$, so that $\R^k$ is regular.

Let $Z\subset (T^* M\oplus\Lambda^2 T^* M)\otimes\R^{n-h}$ denote the set of triples $(m,\alpha,\omega)$ with $\alpha\in T_{m}^{*}M\otimes\R^{n-h}$ and $\omega\in \Lambda^2 T_{m}^{*}M\otimes\R^{n-h}$ such that either $\alpha\in Hom(T_m M,\R^{n-h})$ is not surjective or $\omega_{|ker(\alpha)}$ does not admit any regular isotropic $k$-subspaces. Then $Z$ is a finite union of proper submanifolds.

Let $\Omega$ denote the space of smooth $\R^{n-h}$-valued differential 1-forms $\theta$ on $M$. The map \begin{eqnarray*}
\Omega\times M\to (T^* M\oplus\Lambda^2 T^* M)\otimes\R^{n-h},\quad (\theta,m)\mapsto (\theta(m),d\theta(m))
\end{eqnarray*}
is transverse to $Z$. In fact, in case $M=\R^n$, the restriction of this (linear) map to the finite dimensional space of differential forms of the form
\begin{eqnarray*}
\theta_{a,b,y}=\sum_{i}a_{i}dx_i +\sum_{i,j}b_{i,j}(x_i -y_i)dx_j
\end{eqnarray*}
is already onto at each point, i.e. a submersion. As a consequence, for a generic choice of $\theta$, the section $(\theta,d\theta)$ is transverse to $Z$. The corresponding $H=ker(\theta)$ admits $k$-dimensional regular isotropic subspaces at each point except those of $(\theta,d\theta)^{-1}(Z)$, a union of proper submanifolds, which has measure zero.\qed 

\subsection{Infinitesimal existence of regular horizontal immersions}

If a smooth germ of immersion $f:(\R^k ,0)\to (M,m)$ satisfies the horizontality equation up to order 1, i.e. $E(f)(v)=o(|v|^1 )$, then $S=im(d_0 f)$ is isotropic. If $S$ turns out to be regular, there is no other algebraic obstruction to deforming $f$ to a horizontal immersion, at least at the level of finite jets.

\begin{prop}
\label{infinitesimal}
Let $(M,H)$ be a Carnot manifold. Let $m\in M$ and let $S\subset H_m$ be a $k$-dimensional regular isotropic subspace. Then there exists a germ of immersion $f:(\R^k ,0)\to (M,m)$ which satisfies the horizontality equation to infinite order, i.e.
\begin{eqnarray*}
E(f)(v)=o(|v|^N ) \quad \textrm{for all integers}\quad N,
\end{eqnarray*}
and such that $im(df)=S$ at the origin.
\end{prop}

\proof
Choose coordinates on $V=\R^k$ and $M=\R^n$. Let us prove, by induction on $N$, that there exists a polynomial $f_N$ of degree $N+1$ such that $f(0)=0$, $im(d_0 f)=S$ and
\begin{eqnarray*}
f_{N}^{*}\theta(v)=o(|v|^N ),\quad f_{N}^{*}d\theta(v)=o(|v|^N ).
\end{eqnarray*}

Any linear immersion $f_0$ such that $im(f_0 )=S$ satisfies 
\begin{eqnarray*}
f_{0}^{*}\theta(v)=o(|v|^0 ),\quad f_{0}^{*}d\theta(v) =o(|v|^0 ).
\end{eqnarray*}

Let $N\geq 1$. Assume $f_{N-1}$ exists. Since  for arbitrary $f$, $f^*\theta(\frac{\partial f}{\partial x_i})=\theta(f)(\frac{\partial f}{\partial x_i})$ is a product, when differentiating $f^*\theta$ $N$ times, all terms but one involve less than $N+1$ derivatives of $f$. Therefore, if $f=f_{N-1}+o(|v|^{N+1})$, then, at the origin,
\begin{eqnarray*}
\frac{\partial^{N}}{\partial x_{i_{1}}\cdots\partial x_{i_{N}}}(f^*\theta(\frac{\partial}{\partial x_{i_{N+1}}}))=\theta(f)  (\frac{\partial^{N+1}f}{\partial x_{i_{1}}\cdots\partial x_{i_{N+1}}})+\frac{\partial^{N}}{\partial x_{i_{1}}\cdots\partial x_{i_{N}}}(f_{N-1}^* \theta(\frac{\partial}{\partial x_{i_{N+1}}})),
\end{eqnarray*}
and
\begin{eqnarray*}
&&\frac{\partial^{N}}{\partial x_{i_{1}}\cdots\partial x_{i_{N}}}(f^*d\theta(\frac{\partial}{\partial x_{i_{N+1}}},\frac{\partial}{\partial x_{i_{N+2}}}))
=(d\theta)(f) (\frac{\partial^{N+1}f}{\partial x_{i_{1}}\cdots\partial x_{i_{N+1}}},\frac{\partial f}{\partial x_{i_{N+2}}})\\
&+&
d\theta(f)(\frac{\partial f}{\partial x_{i_{N+1}}},\frac{\partial^{N+1}f}{\partial x_{i_{1}}\cdots\partial x_{i_{N}}\partial x_{i_{N+2}}})
+\frac{\partial^{N}}{\partial x_{i_{1}}\cdots\partial x_{i_{N}}}(f_{N-1}^*d\theta(\frac{\partial}{\partial x_{i_{N+1}}},\frac{\partial}{\partial x_{i_{N+2}}})).
\end{eqnarray*}
We seek $f_N$ in the form
\begin{eqnarray*}
f_N (v)=f_{N-1}(v) +X(v,\ldots,v),
\end{eqnarray*}
where $X\in S^{N+1}V^* \otimes TM$ is an unknown $TM=\R^n$-valued symmetric tensor. The equations to be solved are of the form
\begin{eqnarray*}
\theta \circ X=-a,\quad \mathcal{A}(\alpha\circ X)=-b,
\end{eqnarray*}
with the following notation. 
\begin{itemize}
  \item $\theta$, evaluated at the origin, belongs to $Hom(TM,\R^{n-h})$, thus can be composed with $X$ to yield $\theta\circ X\in S^{N+1}V^* \otimes\R^{n-h}$.
  \item $a=\sum_{i_1 ,\ldots,i_{N+1}}\frac{\partial^{N}}{\partial x_{i_{1}}\cdots\partial x_{i_{N}}}(f_{N-1}^*\theta(\frac{\partial}{\partial x_{i_{N+1}}}))dx_{i_1}\cdots dx_{i_N}\otimes dx_{N+1}\in S^{N}V^* \otimes V^* \otimes\R^{n-h}$.
  \item $\alpha\in Hom(TM,V^* \otimes\R^{n-h})$ is defined by $\alpha(w)=\iota_w d\theta\circ f_0$. It can be composed with $X$ to yield $\alpha \circ X\in S^{N+1}V^* \otimes V^* \otimes\R^{n-h}$.
  \item $\mathcal{A}:S^{N+1}V^* \otimes V^* \otimes \R^{n-h}\to S^{N}V^* \otimes \Lambda^2 V^* \otimes \R^{n-h}$ denotes skew-symmetrization with respect to the last two variables. It maps an $\R^{n-h}$-valued tensor $T$ to $\mathcal{A}(T):(v_1 ,\ldots,v_{N+2})\mapsto -T(v_1 ,\ldots,v_{N+1},v_{N+2})+T(v_1 ,\ldots,v_N ,v_{N+2},v_{N+1})$.
  \item $b=\sum_{i_1 ,\ldots,i_{N+2}}\frac{\partial^{N}}{\partial x_{i_{1}}\cdots\partial x_{i_{N}}}(f_{N-1}^*d\theta(\frac{\partial}{\partial x_{i_{N+1}}},\frac{\partial}{\partial x_{i_{N+2}}}))dx_{i_1}\cdots dx_{i_N} \otimes dx_{N+1}\wedge dx_{N+2}\in S^{N}V^* \otimes \Lambda^2 V^*\otimes\R^{n-h}$.
\end{itemize}
For these linear equations to admit solutions, there are two necessary conditions : $a$ should be fully symmetric (for this, $\mathcal{A}(a)=0$ suffices) and $b$ should satisfy $\mathcal{C}(b)=0$, where 
$\mathcal{C}(b)\in S^{N-1}V^* \otimes \Lambda^3 V^*\otimes\R^{n-h}$ is given by
\begin{eqnarray*}
\mathcal{C}(b)(v_1 ,\ldots,v_{N+2})
&=&b(v_1 ,\ldots,v_{N+2})+b(v_1 ,\ldots,v_{N-1},v_{N+1},v_{N+2},v_N )\\
&&+b(v_1 ,\ldots,v_{N-1},v_{N+2},v_N ,v_{N+1}).
\end{eqnarray*}
Indeed, $\mathcal{C}\circ\mathcal{A}=0$. 

By definition, $a(v_1 ,\ldots,v_{N+1})=v_1 \cdots v_{N}(f_{N-1}^* \theta(v_{N+1}))$,
\begin{eqnarray*}
\mathcal{A}(a)(v_1 ,\ldots,v_{N+1})&=&
a(v_1 ,\ldots,v_{N-1},v_{N},v_{N+1})-a(v_1 ,\ldots,v_{N-1},v_{N+1},v_N )\\
&=&
v_1 \cdots v_{N-1}(v_N (f_{N-1}^* \theta(v_{N+1}))-v_{N+1}(f_{N-1}^* \theta(v_{N})))\\
&=&v_1 \cdots v_{N-1}(f_{N-1}^* d\theta(v_N ,v_{N+1}))\\
&=&0,
\end{eqnarray*}
since $f_{N-1}^* d\theta=o(|v|^N )$.

In the same way, $b(v_1 ,\ldots,v_{N+2})=v_1 \cdots v_{N}(f_{N-1}^* d\theta(v_{N+1},v_{N+2}))$,
\begin{eqnarray*}
\mathcal{C}(b)(v_1 ,\ldots,v_{N+2})&=&
v_1 \cdots v_{N-1}(v_{N}(f_{N-1}^* d\theta(v_{N+1},v_{N+2}))\\
&&+ v_{N+1}(f_{N-1}^* d\theta(v_{N+2},v_{N}))
+v_{N+2}(f_{N-1}^* d\theta(v_{N},v_{N+1})))\\
&=&v_1 \cdots v_{N-1}(f_{N-1}^* dd\theta(v_{N},v_{N+1},v_{N+2}))\\
&=&0,
\end{eqnarray*}
since $dd\theta=0$.

Since $\theta\in Hom(TM,\R^{n-h})$ is onto, it admits a right inverse $\theta^{-1}$, $\theta\circ\theta^{-1}=id_{\R^{n-h}}$. Then $Z=-\theta^{-1}\circ a$ 
satisfies $\theta\circ Z=-a$. In order to solve simultaneously the second equation $\mathcal{A}(\alpha\circ X)=-b$, we look for $Y\in S^{N+1}V^* \otimes \ker(\theta)$ such that 
\begin{eqnarray*}
\mathcal{A}(\alpha\circ Y)=-b-\mathcal{A}(\alpha\circ Z).
\end{eqnarray*}
Thanks to the regularity assumption, the restriction of $\alpha$ to $\ker(\theta)\to V^* \otimes\R^{n-h}$ is surjective. Pick a right inverse $\alpha^{-1}:V^* \otimes\R^{n-h}\to\ker(\theta)$, $\alpha\circ\alpha^{-1}=id_{V^* \otimes\R^{n-h}}$. Look for $Y$ in the form $Y=\alpha^{-1}\circ Y'$ for $Y'\in S^{N+1}V^* \otimes V^* \otimes\R^{n-h}$, which must satisfy
\begin{eqnarray*}
\mathcal{A}(Y')=-b-\mathcal{A}(\alpha\circ Z).
\end{eqnarray*}
Note that $\mathcal{C}(-b-\mathcal{A}(\alpha\circ Z))=0$. According to Lemma \ref{koszul}, this is a sufficient condition for the existence of a solution $Y'$. This completes the inductive proof. 

Once the infinite jet $f_{\infty}$ is found, any germ $(\R^k ,0)\to (M,m)$ having this power series as Taylor expansion satisfies $E(f)=o(|v|^N )$ for all $N$.\qed

\begin{lemme}
\label{koszul}
The sequence 
\begin{eqnarray*}
0\to S^{N+2}V^* \hookrightarrow S^{N+1}V^* \otimes V^* \fleche{\mathcal{A}} S^{N}V^* \otimes \Lambda^2 V^* \fleche{\mathcal{C}} S^{N-1}V^* \otimes \Lambda^3 V^*
\end{eqnarray*}
is exact.
\end{lemme}

\proof
It is a subcomplex of the de Rham complex of $V=\R^k$. Indeed, elements of $\bigoplus_{\ell}S^{N-\ell}V^* \otimes \Lambda^{\ell}V^*$, viewed as differential forms with polynomial coefficients, are exactly the smooth differential forms on $V$ which are homogeneous of degree $N$. Therefore, they form a subcomplex. The Poincar\'e homotopy formula for solving $d$ is homogeneous, thus the subcomplex is acyclic. Finally, up to a factor of 2, $\mathcal{A}$ and $\mathcal{C}$ coincide with exterior differentials.\qed

\subsection{Algebraic inverses}

Passing from an infinite power series to a true locally converging solution requires some analysis. We shall use an implicit function theorem. As usual, invertibility of the differential of the equation is needed.

\begin{prop}
\label{inverse}
Let $(M,H)$ be a Carnot manifold. If $f:V\to M$ is a regular horizontal immersion (i.e. $D_v f(T_v V)$ is a regular subspace of $H_{f(v)}$ for all $v\in V$), then $D_f E$ admits an algebraic right inverse.
\end{prop}

\proof
It suffices to right invert the map $X\mapsto f^* (\iota_X d\theta)$ pointwise. 

To show that such an inverse can be chosen smoothly, consider the set $RegIso$ of pairs $(m,S)$, $m\in M$, $S$ a regular isotropic subspace of $H_m$. This is a submanifold in the bundle of Grassmannians $Gr(k,H)$. This implies that the set $RegIsoImm$ of triples $(v,m,L)$ where $v\in V$, $m\in M$ and $L:T_v V \to H_m$ is an injective linear map with regular isotropic image is a submanifold in the bundle $Hom(TV,TM)$ over $V\times M$. For $(v,m,L)\in RegIsoImm$, the set $Right_{(v,m,L)}$ of right inverses of the map
\begin{eqnarray*}
H_{m}\to Hom(T_v V,\R^{n-h})=T_{v}^{*}V\otimes \R^{n-h},\quad X\mapsto L^* (\iota_X d\theta),
\end{eqnarray*}
is an affine space of constant dimension $(n-h)k(h-(n-h)k)$. The spaces $Right_{(v,m,L)}$ form a smooth bundle with contractible fibers, therefore it admits a smooth section $(v,m,L)\mapsto right(v,m,L)$ defined on $RegIso$. 

For $g$ a $\R^{n-h}$-valued 1-form on V, consider the vector field along $f$ defined by $X(v)=right(v,f(v),d_v f)(g(v))$. The map $M_f :g\mapsto X$ is a right inverse of $D_f E$.\qed

\begin{rem}
\label{reminverse}
The existence of an algebraic right inverse for a differential operator is not unusual. In fact, it is generic for \emph{underdetermined} operators, i.e., with more unknown functions than equations, see \cite{Gromov-PDR}, page 156. 
\end{rem}

\begin{exemple}
\label{exinverse}
\emph{(J. Nash).} A specific example is the linearization at a \emph{free} map of the isometric immersion operator which to a map $f$ between Riemannian manifolds associates the pulled back metric. 
\end{exemple}
A map $\R^k \to \R^n$ is \emph{free} if all its first and second derivatives are linearly independant. This notion was introduced by E. Cartan and M. Janet, see \cite{Janet}.

\begin{rem}
\label{leftinverse}
Generic \emph{overdetermined} operators admit differential left inverses, see \cite{Gromov-PDR}, page 166.
\end{rem}

\section{Implicit function theorem}

\subsection{Nash's implicit function theorem}

In his solution of the isometric embedding problem, J. Nash discovered that differential operators whose linearization admits a differential right inverse can be right inverted. This implicit function theorem has an unusual feature: the inverse map is a \emph{local} operator.

\begin{theo}
\label{nash}
\emph{(J. Nash. This version is taken from \cite{Gromov-PDR}, page 117).} Let $V$ be a Riemannian manifold. Let $F$, $G$ be bundles over $V$. Let $E:C^{\infty}(F)\to C^{\infty}(G)$ be a differential operator whose linearization $D_f E$ admits a differential right inverse $M_f$, which is defined for $f$ in a subset $\mathcal{A}$ of $C^{\infty}(F)$ defined by strict differential inequalities. Fix a real number $\rho>0$. Then
there exist an integer $s$ such that the following holds.

For each $f\in \mathcal{A}$, there exists a right inverse $E_{f}^{-1}$ of $E$, defined on a $C^s$-neighborhood of $E(f)$ in $C^{\infty}(G)$, such that $E_{f}^{-1}(E(f))=f$ and whose differential at $E(f)$ is $M_f$. Furthermore, $E_{f}^{-1}$ depends smoothly on parameters, and is \emph{local}: given $f\in \mathcal{A}$, $g\in C^{\infty}(G)$ and $v\in V$, $E_{f}^{-1}(g)(v)$ depends only on the values of $f$ and $g$ in a ball of radius $\rho$ centered at $v$.
\end{theo}

\proof
It can be found in textbooks like \cite{Gromov-PDR} or \cite{AG}.\qed

\subsection{Local existence}

\begin{cor}
\label{approxsol}
\emph{(\cite{Gromov-PDR} page 119).} Same assumptions as in Theorem \ref{nash}. Any germ $f_0$ that solves
\begin{eqnarray*}
E(f_0 )(v)=o(|v-v_0 |^s )
\end{eqnarray*}
can be deformed to a true local solution $f_1$: $E(f_1 )=0$.
\end{cor}

Indeed, choose $g\in C^{\infty}(G)$ such that $g=-E(f_0 )$ near $v_0$, but $g$ is $C^s$-small. For $t\in[0,1]$, set $f_t =E_{f_0}^{-1}(E(f_0 )+tg)$.\qed

In other words, it suffices to construct solutions up to order $s$ ($s=2$ is enough for the horizontal immersion problem). With Proposition \ref{infinitesimal}, this completes the proof of the existence of local regular horizontal immersions. It even gives a more precise information on the topology of the space of germs of solutions at a point. 

\begin{prop}
\label{locex}
Let $(M,H)$ be a Carnot manifold. Let $m\in M$. Consider the space $Sol_m$ of germs of regular horizontal immersions $(\R^{k},0)\to (M,m)$. Map a germ to its derivative at the origin. This gives a homotopy equivalence of $Sol_m$ to the space $InjRegIso_m$ of injective linear maps $\R^k \to H_m$ with regular isotropic image.
\end{prop}

\proof
The proof of Proposition \ref{infinitesimal} shows that to construct a solution to order $s$ from an injective linear map $\R^k \to H_m$ with regular isotropic image, one merely needs solve linear equations, i.e. pick points in nonempty affine spaces of constant dimension. Therefore, the space $Jet_m=\{\textrm{germs}\,f\,|\, f(0)=m,\, E(f)=o(|v|^s )\}$ is homotopy equivalent to $InjRegIso_m$. The map $[0,1]\times Jet_m \to Jet_m$, $(f_0 ,t)\mapsto (f_t )$ provided by Corollary \ref{approxsol} is a deformation retraction of $Jet_m$ to $Sol_m$.\qed  

\subsection{Microflexibility}

Here is a second application of Nash's implicit function theorem. When it applies, solutions can evolve more or less independantly on disjoint parts of their domains. This vague statement is made precise in the following definition.

\begin{defi}
\label{flex}
\emph{(\cite{Gromov-PDR} page 41).} Say an equation is \emph{flexible} (resp. \emph{microflexible}) if, given compact sets $K'\subset K\subset V$, a solution $f$ defined on a neighborhood of $K$ and a deformation $f_t$, $t\in[0,1]$, of its restriction to a neighborhood of $K'$, the deformation extends to a neighborhood of $K$ (resp. for a little while, i.e. for $t\in[0,\epsilon]$). We also require a parametrized version of this property: it should apply to continuous families of solutions $f_p$ and of deformations $f_{p,t}$ parametrized by an arbitrary polyhedron $P$.
\end{defi}

\begin{exemple}
\label{exflex}
Strict inequations are trivially microflexible, but need not be flexible.
\end{exemple}
Apart from this trivial example, microflexibility is hard to establish. The following result gives all examples I know.

\begin{cor}
\label{inversemicro}
\emph{(\cite{Gromov-PDR} page 120).} Let $E$ be a differential operator such that $D_f E$ admits a differential right inverse for $f\in \mathcal{A}$. Then the system $\mathcal{A}\cap\{E=0\}$ is microflexible.
\end{cor}

\proof
Given solutions $f$ near $K$ and $f_t$ near $K'$, construct a family of sections $f'_t$ defined near $K$ which coincides with $f_t$ on the $\rho$-neighborhood of $K'$ for some $\rho>0$. For $t$ small, one can set $e_t =E_{f'_t}^{-1}(0)$. It is a solution defined near $K$. Furthermore, near $K'$, $e_t =f_t$ by locality.\qed

\begin{exemple}
\label{exinversemicro}
In the underdetermined case $h-k\geq (n-h)k$, for a generic distribution (resp. for a contact distribution), $k$-dimensional horizontal immersions are microflexible.
\end{exemple}

\begin{rem}
\label{remfibration}
Flexibility means that the restriction map $f\mapsto f_{|neigh(K')}$ between spaces of solutions near $K$ (resp. near $K'$) is a \emph{Serre fibration} (path lifting property).

Microflexibility sounds like this map being a submersion. 
\end{rem}
As we shall see in the next section, this opens the way to topological methods for the study of homotopy properties of spaces of solutions.

\subsection{Calculus of variations}

Before we proceed, let us mention our last consequence of Nash's implicit function theorem. 

\begin{lemme}
\label{calcul}
\emph{(\cite{Gromov-CC} page 254).} Same assumptions as in Theorem \ref{nash}. There is an open neighborhood $\mathcal{U}$ of $\mathcal{A} \cap\{E=0\}$ in $\mathcal{A}$ and a smooth retraction $r:\mathcal{U}\to\mathcal{A}\cap\{E=0\}$ with differential
\begin{eqnarray*}
D_f r=1-M_f \circ D_f E .
\end{eqnarray*}
\end{lemme}

\proof
The domain $\mathcal{V}$ of $E^{-1}:(f,g)\mapsto E_{f}^{-1}(g)$ is an open $C^s$-neighborhood of $\{(f,g)\,|\,E(f)=g\}$ in $\mathcal{A}\times C^{\infty}(G)$. Let $\mathcal{U}=\{f\in\mathcal{A}\,|\,(f,0)\in\mathcal{V}\}$. For $f\in\mathcal{U}$, let $r(f)=E^{-1}(f,0)$. Then $E(r(f))=0$. If $f\in\mathcal{A}$ and $E(f)=0$, then $f\in\mathcal{U}$ and $r(f)=f$.

By construction, $\frac{\partial E^{-1}}{\partial g}(f,E(f))=M_f$. Since $E^{-1}(f,E(f))=f$ for all $f\in\mathcal{A}$,
\begin{eqnarray*}
\frac{\partial E^{-1}}{\partial f}(f,E(f))+\frac{\partial E^{-1}}{\partial g}(f,E(f))\circ D_f E=id,
\end{eqnarray*}
and the formula for $D_f r$ follows.\qed

\medskip

It follows that given a smooth functional $\Phi$ on $\mathcal{A}$, one can assert that
\begin{eqnarray*}
f \textrm{ critical point of }\Phi_{|\mathcal{A}\cap\{E=0\}}\quad\Leftrightarrow\quad D_f \Phi \textrm{ vanishes on }im(D_{f}r),
\end{eqnarray*}  
and derive Euler-Lagrange equations for the restriction of $\Phi$ to $\mathcal{A}\cap\{E=0\}$ in the usual manner.

\medskip

This applies to the area functional for \emph{regular} horizontal immersions $f$ in a Carnot manifold. With the algebraic inverse $M_f$ constructed in Proposition \ref{inverse}, $D_f r$ maps a vector field $X$ along $f$ to an $X'=X+Y$ where $Y$ is horizontal and is algebraicly chosen to satisfy $f^* (\iota_{X'}d\theta)=-d(\iota_X \theta)$. In other words, $\delta=D_f r$ is a first order linear differential operator. Since $X$ can be chosen with arbitrarily small support, the usual integration by parts yields the third order equation $\delta^* h=0$ where $\delta^*$ is the adjoint of $\delta$ and $h$ the Riemannian mean curvature of the immersion.

\medskip

The case of contact 5-manifolds is of particular interest, see for example \cite{Schoen-Wolfson}, \cite{Romon}. 

\section{From local to global horizontal immersions}

\subsection{Sketch of proof of Theorem \ref{generich}}
\label{sketch}

Let us sketch a proof of Theorem \ref{generich} for $k=1$. 

Start with a continuous loop $l_0 :P=\R/\Z\to M$. Assume it admits a continuous lift $F_0$ which chooses for each $p\in P$ a regular isotropic line in $H_{f_0 (p)}$.  Proposition \ref{locex} allows a continuous choice of a germ $\bar{f}_p$ of regular horizontal immersion $(\R,0)\to (M,f_0 (p))$, tangent to $F_0 (p)$ at $f_0 (p)$ for each $p\in P$.

\begin{center}
\includegraphics[width=2in]{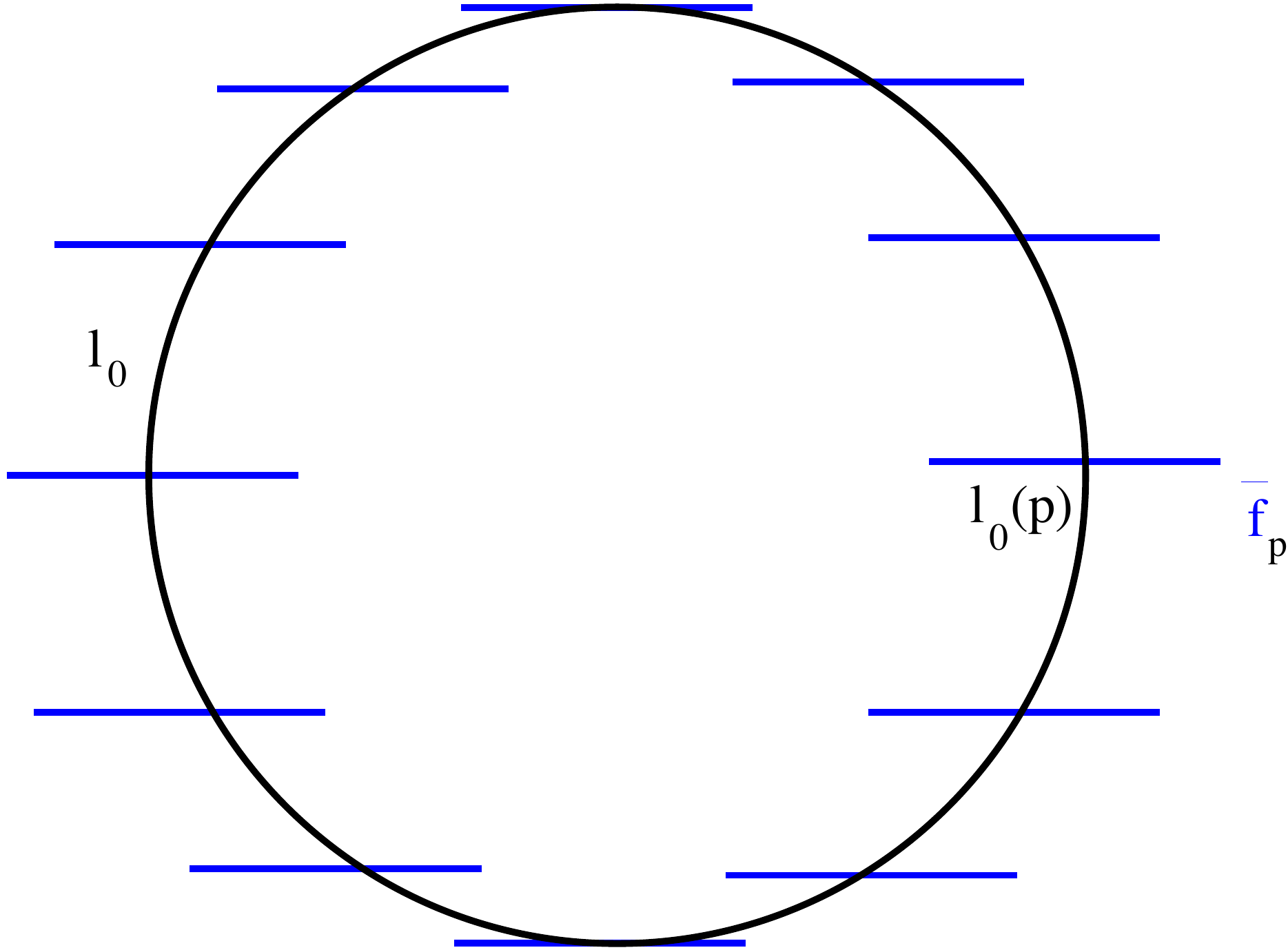}
\par
A curve of germs of regular horizontal curves
\end{center}

Say $\bar{f}_{p}$ is defined on $(p-1,p+1)$. Fix some $p$. Let $K'=K'_+ \cup K'_-$ where $K'_{+}=[p+.5,p+.9]$ and $K'_{+}=[p-.9,p-.5]$, and $K=[p-.9,p+.9]$. Apply microflexibility to the parametric deformation $f_{p,t}$ such that $f_{p,0}=\bar{f}_p$, $f_{p,t} =\bar{f}_{p+t}$ near $K'_+$ and $f_{p,t} =\bar{f}_{p}$ near $K'_-$. 

\begin{center}
\includegraphics[width=2in]{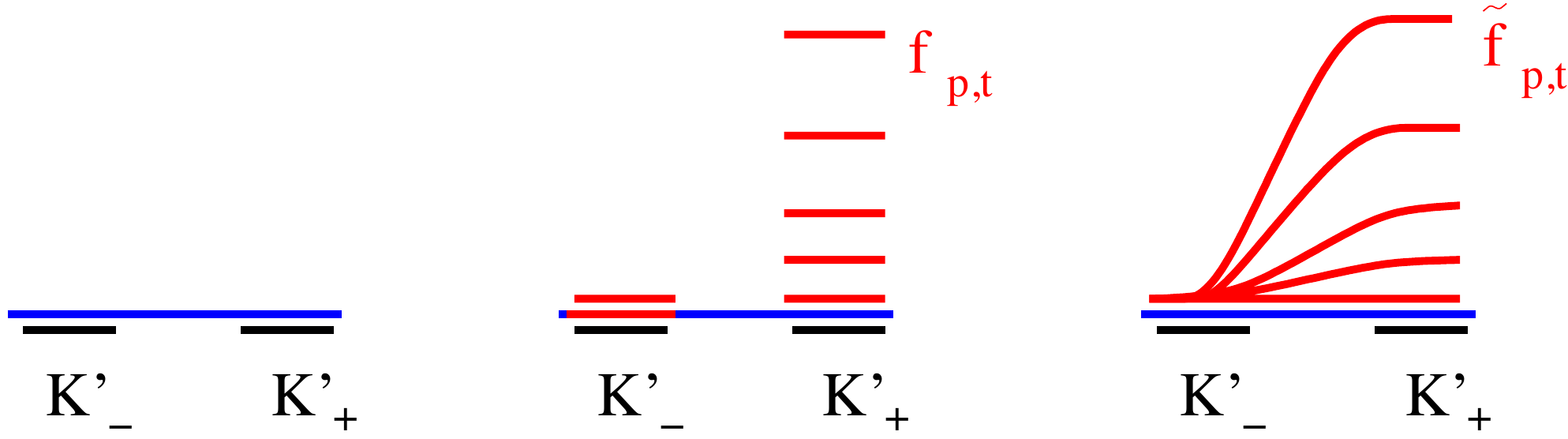}
\par
A deformation of solutions provided by microflexibility
\end{center}

This gives an $\epsilon>0$ and a family of solutions $\tilde{f}_{p,t}$, $t\in[0,\epsilon]$. $\tilde{f}_{p+\epsilon}$ interpolates between $\bar{f}_p$ and $\bar{f}_{p+\epsilon}$ (see figure). Divide $\R/\Z$ in $1/\epsilon$ intervals, glue them together along the pattern shown on the figure. One gets a 1-dimensional branched manifold, together with a regular horizontal immersion to $M$ which is $\epsilon$-close to $l_0$.

\begin{center}
\includegraphics[width=2in]{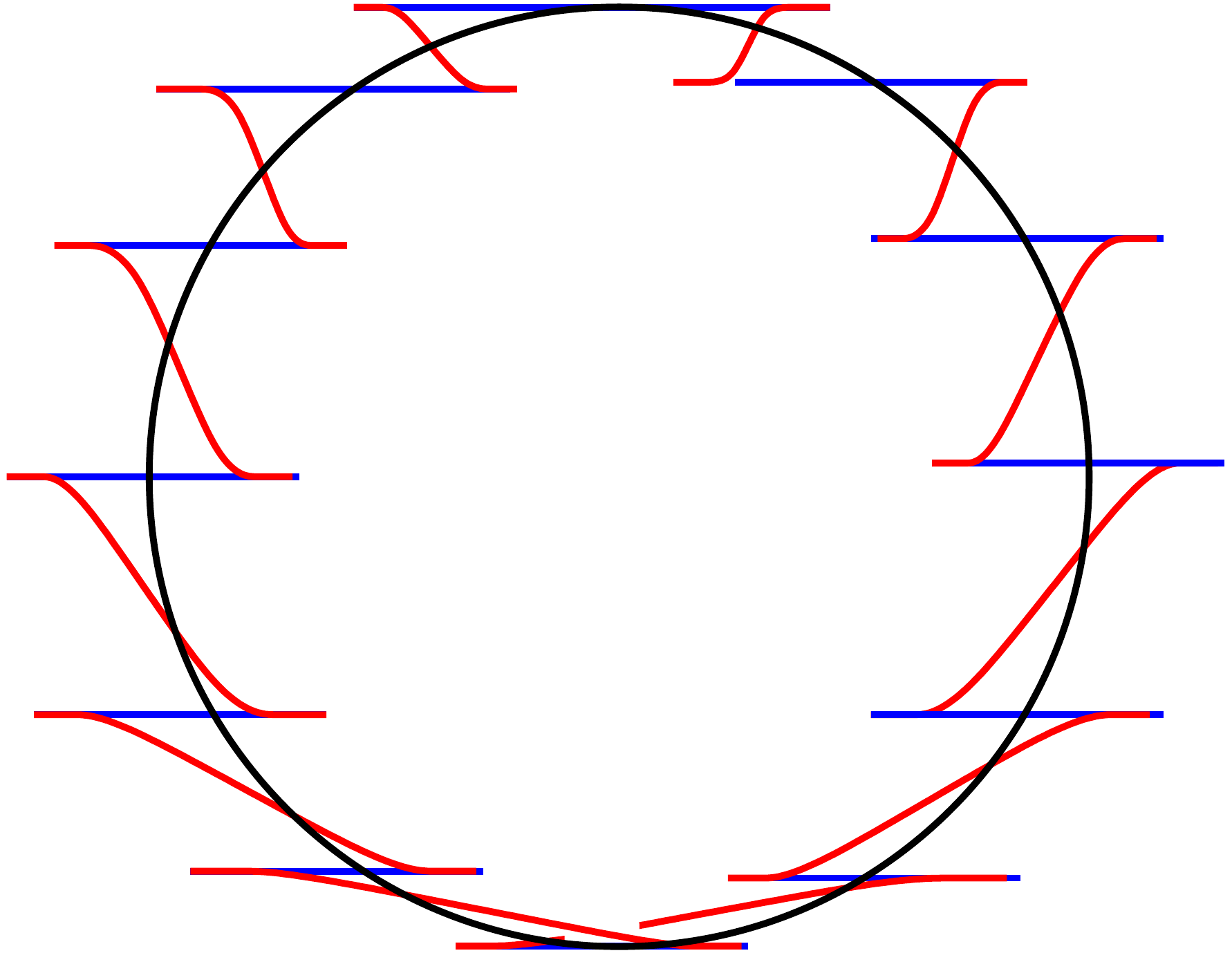}
\par
A branched horizontal immersion
\end{center}

The higher dimensional case requires more topological formalism, some of it is explained below. 

In fact, the technique comes very close to producing global (unbranched) horizontal immersions from honest manifolds. Although this is not needed for our application to dimensions of submanifolds, we shall collect enough material to prove a density property of horizontal immersions when the domain is an open manifold. 

\subsection{Nonholonomic solutions}

\begin{defi}
\label{nonholo}
Let $F\to V$ be a smooth bundle. A differential equation of order $r$ with unknowns in $C^{\infty}(F)$ should be viewed as a subset $\mathcal{R}$ in the bundle of $r$-jets of sections of $F$. Sections of $\mathcal{R}$ are called \emph{nonholonomic solutions} of the equation. 
\end{defi}

\begin{exemple}
\label{exnonholo}
Horizontal immersions $V\to (M,H)$.
\end{exemple}
This is a first order equation. A 1-jet of map $V\to M$ at $v\in V$ is a pair $(m,L)$ where $m\in M$ and $L:T_v V\to T_m M$ is linear. Then $\mathcal{R}$ is the subset of pairs $(m,L)$ where $L$ is injective with its image contained in $H_m$. 

A nonholonomic solution consists of a continuous map $f:V\to M$ and an injective bundle map $F:TV\to H$ over $f$, i.e. $F$ maps injectively $T_v V$ to $H_{f(v)}$.

\subsection{$h$-principle}

We want to study a refinement of the existence problem for solutions of an equation: can one approximate nonholonomic solutions with solutions ? Since solutions form a $C^r$ closed set, one cannot require a $C^r$ approximation, but one can still hope for a $C^0$ approximation. For instance, in the case of horizontal immersions, the data is a bundle map $(f_0 ,F_0 )$, and one wants a horizontal immersion $f_1:V\to M$ such that $f_1$ is $C^0$ close to $f_0$.

It turns out to be very fruitful to investigate simultaneously the homotopy type of the space of solutions, i.e. to require the $r$-jet of the approximating solution to be homotopic to the given nonholonomic solution among nonholonomic solutions.

\begin{defi}
\label{hprinciple}
\emph{(\cite{Gromov-PDR} pages 3, 17, 18).}
Say an equation satisfies the \emph{parametric $C^0$-dense $h$-principle} ($h$-principle, for short) if for every nonholonomic solution, there is a $C^0$-small homotopy to a solution (and also familywise).
\end{defi}

The $h$-principle localizes near a compact subset $K$, and in particular near points. 

\begin{prop}
\label{localh}
Regular horizontal immersions satisfy the $h$-principle near points.
\end{prop}

\proof
This is precisely what Proposition \ref{locex} means.\qed

\subsection{$h$-principle as a homotopy theory}

The $h$-principle has a relative version for a pair $(K,K')$ (\cite{Gromov-PDR} page 39): say that $h$-principle holds for the pair if for every nonholonomic solution defined near $K$, which is a solution near $K'$, there is a $C^0$-small homotopy to a solution defined near $K$, where the homotopy is constant near $K'$ (and again one wants this to hold familywise).

\begin{prop}
\label{k'+kk'=k}
\emph{(\cite{Gromov-PDR} pages 40, 42).}
If the $h$-principle holds near $K'$ and for the pair $(K,K')$, then it holds for $K$.

If the $h$-principle holds near $K_1 \cap K_2$ and for the pairs $(K_1 ,K_1 \cap K_2 )$ and $(K_2 ,K_1 \cap K_2 )$, then it holds for $K_1 \cup K_2$.
\end{prop}

\proof
Start with a nonholonomic solution $F_0$ defined on a neighborhood $U$ of $K$. The $h$-principle near $K'$ provides a $C^0$-small homotopy $F_t$ defined on a neighborhood $U'$ of $K'$, with $F_1 =j^r f_1$ a solution. Extend $F_1$ into a nonholonomic solution $\tilde{F}_1$ defined on $U$ as follows: using a tubular neighborhood of $\partial U'$, write
\begin{eqnarray*}
U=U'\cup(\partial U'\times[0,1])\cup U\setminus U',
\end{eqnarray*}
and set
\begin{eqnarray*}
\tilde{F}_1 &=&F_1 \quad\textrm{on}\quad U',\\
\tilde{F}_1 (u,t)&=&F_{1-t}(u) \quad\textrm{on}\quad \partial U'\times[0,1],\\
\tilde{F}_1 &=&F_0 \quad\textrm{on}\quad U\setminus U'.
\end{eqnarray*}
The relative $h$-principle finally provides a $C^0$-small homotopy of $\tilde{F}_1$ to a solution defined near $K$. This works for families as well, thus the $h$-principle holds near $K$.

In case $K'=K_1 \cap K_2$, the same construction gives a $C^0$-small homotopy to a solution defined near $K_1 \cup K_2$, thanks to the extension character of the relative $h$-principle.\qed

\begin{theo}
\label{k'+k=kk'}
\emph{(S. Smale. Taken from \cite{Gromov-PDR} page 42).} Assume that 
\begin{itemize}
  \item the $h$-principle holds for $K$ and $K'$,
  \item the equation is flexible.
\end{itemize}
 Then the $h$-principle holds for the pair $(K,K')$.
\end{theo}

\proof
Consider the diagram
\begin{eqnarray*}
\begin{array}{ccc}
\{\textrm{solutions near K}\}&\to&\{\textrm{nonholonomic solutions near K}\}\\
\downarrow&&\downarrow \\
\{\textrm{solutions near K'}\}&\to&\{\textrm{nonholonomic solutions near K'}\}
\end{array},
\end{eqnarray*}
where horizontal arrows are inclusions and vertical arrows restriction maps. Since the $h$-principle holds near $K$ and near $K'$, horizontal arrows are weak homotopy equivalences. Flexibility makes vertical arrows Serre fibrations. From the five lemma applied to the long homotopy exact sequences of these fibrations, we know that the fibers are weakly homotopic. In particular, if a solution near $K'$ admits an extension as a noholonomic solution, it also extends as a solution, and this works for families as well. This is the relative $h$-principle.\qed

\begin{cor}
\label{flexh}
\emph{(\cite{Gromov-PDR} page 42).} Assume that 
\begin{itemize}
  \item the $h$-principle holds near points,
  \item the equation is invariant under diffeomorphisms of the domain,
  \item the equation is flexible.
\end{itemize}
 Then the $h$-principle holds in $V$.
 \end{cor}

\proof
Triangulate $V$. By assumption, the $h$-principle holds in balls embedded in $V$. Since each simplex has a basis of neighborhoods diffeomorphic to a ball, the $h$-principle holds near simplices. Flexibility and Theorem \ref{k'+k=kk'} imply that $h$-principle holds for all pairs. Proposition \ref{k'+kk'=k} implies that it passes to unions, i.e. to all of $V$.\qed

\subsection{From microflexibility to flexibility}

There remains to establish flexibility for certain equations. It turns out that, for equations which are invariant on reparametrization of the domain, microflexibility implies flexibility in one dimension less.

\begin{theo}
\label{micro=>flex}
\emph{(\cite{Gromov-PDR} page 78).}
Consider an equation on $V$ which is $Diff(V)$-invariant and microflexible. Then flexibility holds for germs of solutions along any proper submanifold.
\end{theo}

The proof is sketched in the next paragraph.

\begin{cor}
\label{open}
\emph{(\cite{Gromov-PDR} page 79).}
A microflexible $Diff(V)$-invariant equation which satisfies the $h$-principle near points satisfies the $h$-principle on open manifolds.
\end{cor}

Indeed, if $V$ is open, there exists a codimension 1 polyhedron $V_0$ in $V$ and an isotopy which maps every compact subset of $V$ into arbitrarily small neighborhoods of $V_0$. According to Theorem \ref{micro=>flex}, microflexibility on $V$ implies flexibility for solutions defined near $V_0$, which implies flexibility for all solutions. One concludes with corollary \ref{flexh}.\qed

\subsection{Compressibility}

For the proof of theorem \ref{micro=>flex}, it is convenient to replace flexibility by the equivalent notion on compressibility.

\begin{defi}
\label{defcompr}
\emph{(\cite{Gromov-PDR} page 80).} 
Given an equation on a manifold $V$ and a compact set $K\subset V$, call \emph{deformation on $K$} any curve of solutions $f_t$, $t\in[a,b]$, defined in a neighborhood of $K$. Say a deformation $f_t$ on $K$ is \emph{compressible} if for every sufficiently small neighborhood $U$ of $K$, it can be extended to a deformation $\tilde{f}_t$ defined on $U$ for $t\in[a,b]$, which is constant (i.e. independant on $t$) in a neighborhood of $\partial U$. One also needs a parametric version of this definition, for parametric families of solutions $f_{t,p}$, $p\in P$.
\end{defi}

\begin{center}
\includegraphics[width=3in]{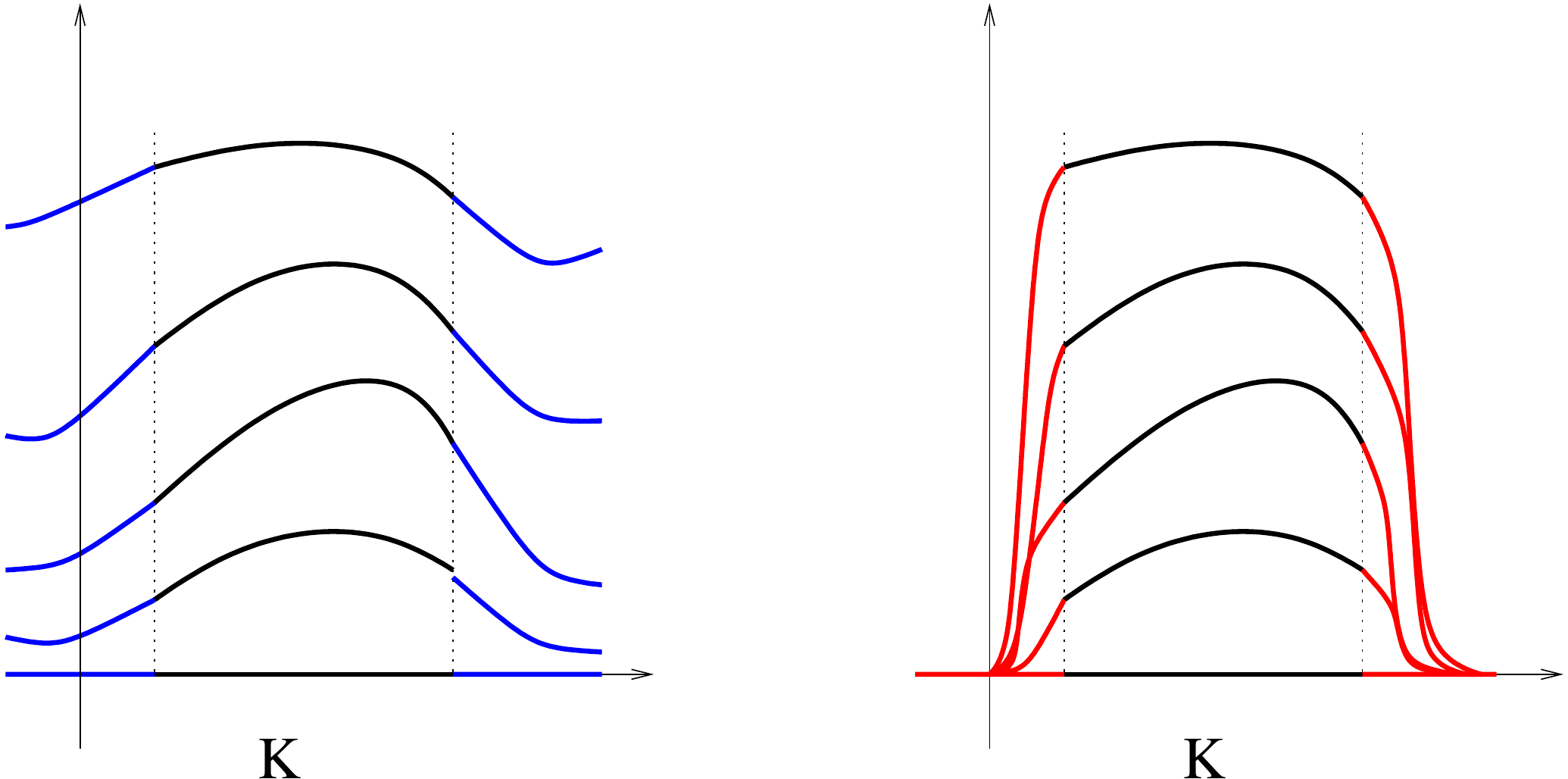}
\par
A compressible deformation
\end{center}

\begin{lemme}
\label{compr<=>flex}
An equation is flexible if and only if all deformations on all compact sets are compressible.
\end{lemme}

\proof
Given $K'\subset K$ and a solution near $K$, a compressible deformation of it on $K'$ trivially extends to $K$ by making it constant outside a neighborhood of $K'$. Conversely, given a deformation on $K$ and a compact neighborhood $U$ of $K$, apply flexibility to the pair $K\subset K\cup \partial U$ to obtain an extension which is constant in a neighborhood of $\partial U$.\qed

\medskip

In the same way, microflexibility implies that every deformation on $K$ can be compressed into an arbitrary neighborhood $U$ of for some time $\epsilon$ depending on the deformation and on $U$. For equations that are invariant under reparametrizations of the domain, and when enough space is available, one can arrange that $\epsilon$ depends on the deformation but not on $U$.

\begin{lemme}
\label{micro=>compr}
\emph{(\cite{Gromov-PDR} page 82).}
Consider an equation on $V$ which is $Diff(V)$-invariant and microflexible. Let $V_0 \subset V$ be a submanifold of positive codimension. Let $K\subset V_0$ be compact. Let $f_t$ be a deformation on $K$. There exists $\epsilon=\epsilon(\{f_t\})$ such that for every neighborhood $U_0$ of $K$ in $V_0$, there exists an extension $\tilde{f}_t$ defined on $U_0$ for $t\in[0,\epsilon]$, which is constant near $\partial U_0$.
\end{lemme}

\proof
Microflexibility applied to $K\subset K\cup\partial U$ yields an extension $\bar{f}_t$ defined on $U$ for $t \in[0,\epsilon]$, $\epsilon=\epsilon(\{f_t\})$, which is constant near $\partial U$. Here, $U$ is a neighborhood of $K$ in $V$.

Given a compact neighborhood $U_0$ of $K$ in $V_0$, one can assume that $\bar{f}_t$ is constant on $V_0 \setminus U_0$ for $t\in[0,\epsilon']$, for some $\epsilon'\leq \epsilon$ which unfortunately depends on $U_0$.

Choose a compact neighborhood $L$ of $U_0$ in $U\cap V_0$, with smooth boundary. There exists an isotopy $\delta_t$ of $V$, $t\in[0,\epsilon]$, which is constant on a neighborhood of $\partial L$, such that $\delta_{\epsilon}$ moves $\partial L$ into the neighborhood of $\partial U$ where $\bar{f}_t$ is always constant. For $t\in[0,\epsilon]$, let 
\begin{eqnarray*}
\tilde{f}_{t}&=&f_{t}\circ\delta_t \quad\textrm{on}\quad L,\\
\tilde{f}_{t}&=&f_{\min\{t,\epsilon'\}}\circ\delta_t
\quad\textrm{on}\quad V_0 \setminus L.
\end{eqnarray*}

\begin{center}
\includegraphics[width=4in]{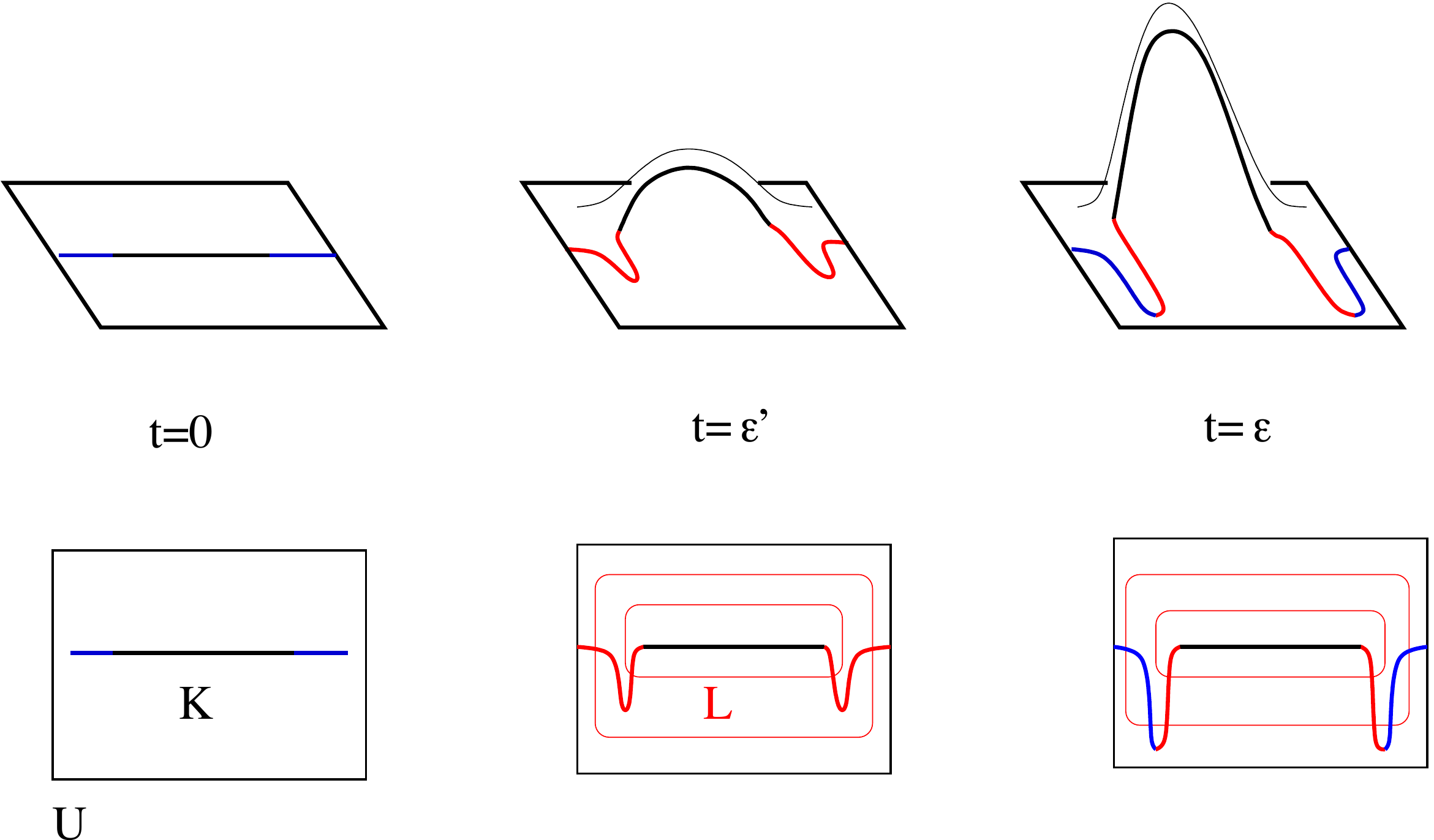}
\par
Exploiting an extra dimension to improve compression
\end{center}

This $\tilde{f}_{t}$ is smooth along $\partial L$. By $Diff(V)$-invariance, $\tilde{f}_{t}$ is a solution. This works more generally for parametric deformations $f_{t,p}$, $p\in P$.\qed

\subsection{Proof of theorem \ref{micro=>flex}}

We must show that every deformation $f_t$, $t\in[0,1]$, on a compact set $K$ is compressible. 

Note that all restrictions of $\{f_t \}$ to smaller intervals $[p,1]$ belong to a unique deformation, the parametric deformation $f_{t,p}$, $p\in P=[0,1]$, defined by
\begin{eqnarray*}
f_{t,p}=f_{\min\{1,t+p\}}.
\end{eqnarray*}
Apply Lemma \ref{micro=>compr} to this parametric deformation. This provides a compression time $\epsilon$ for all deformations $\{f_t \,|\,t\in[p,1]\}$ which depends neither on $p$ nor on a compression neighborhood. 

Fix a neighborhood $U$ of $K$. Compress $\{f_t \,|\,t\in[0,\epsilon]\}$ to $\{\tilde{f}_t \,|\,t\in[0,\epsilon]\}$ within $U$. There exists a neighborhood $U'\subset U$ on which $\tilde{f}_{t}$ coincides with $f_t$ for $t\in[0,\epsilon]$. Compress $\{f_t \,|\,t\in[\epsilon,2\epsilon]\}$ to $\{\tilde{f}_t \,|\,t\in[\epsilon,2\epsilon]\}$ within $U'$. This gives a compression of $\{f_t \,|\,t\in[0,2\epsilon]\}$ within $U'$. Continue.\qed

\subsection{Proof of theorem \ref{generich}}

According to Proposition \ref{isoreggeneric}, in case the polyhedron $P$ is an open manifold, Theorem \ref{generich} follows from the following statement.

\begin{prop}
\label{regiso=>rich}
Let $(M,H)$ be a Carnot manifold. Assume that at some point $m\in M$, there exists a regular isotropic linear subspace $S\subset H_m$ with $dim(S)=k$. Let $f_0$ be continuous map of an open $k$-dimensional manifold $W$ to a sufficiently small neighborhood of $m$. Let $\epsilon>0$. Then there exists a local diffeomorphism $f:W \times \R^{n-k}\to M$ such that all $f_{|W \times\{*\}}$ are horizontal and uniformly $\epsilon$-close to $f_0$.
\end{prop}

\proof
Consider the equation $E$ whose solutions are foliated regular horizontal immersions, i.e. maps $f:V=W \times\R^{n-k}\to M$ such that $d_{(w,z)} f$ maps $T_{w}W$ to a regular isotropic subspace in $H_{f(w,z)}$.

The set $BijRegIso$ of bijective linear maps of $\R^k \times \R^{n-k}$ to $TM$ mapping $\R^k$ to a regular isotropic subspace in some $H_m$ is a submanifold, with $(m',S)\to m'$ a submersion. In some neighborhood $U$ of $m$, one can choose such a bijection $m'\to L(m')$ depending smoothly on $m'$. Then $F(w,z)=(f_0 (w),L\circ f_0 (w))$ is a nonholonomic solution of $E$.

Corollary \ref{inversemicro} implies that $E$ is microflexible. Corollary \ref{open} does not apply directly since the equation is not fully $Diff(V)$-invariant. Nevertheless, all we need is a germ of a solution along the leaf $W\times\{0\}$, so $Diff(W)$-invariance suffices. This allows to $C^0$ approximate $F$ with a foliated horizontal immersion.\qed

\medskip

If $P$ is a closed $n$-manifold, then the same procedure yields $h$-principle along the $n-1$-skeleton of some triangulation. View each $n$-simplex as a 1-parameter family of $n-1$-spheres, and apply the strategy described in paragraph \ref{sketch}, to produce a global solution branched along finitely many spheres.

In fact, the $h$-principle is not needed for the construction of branched solutions, microflexibility is sufficient. We refer to \cite{Gromov-CC} page 262 and \cite{Gromov-PDR}, page 112.\qed

\section{From submanifolds to differential forms}

\subsection{Horizontal forms}

Let $(M,H)$ be a Carnot manifold. Let $p:M\to \R^q$ be a submersion with horizontal fibers. Let $vol$ denote some volume form on $\R^q$. Then the differential form $\eta=p^* vol$ has the following property : if $\theta$ is a 1-form that vanishes on $H$, then $\theta\wedge\eta=0$. This suggests the following definition.

\begin{defi}
\label{defhorform}
Let $\Omega^*$ denote the space of differential forms on $M$. Let $\Theta^* \subset\Omega^*$ denote the ideal of differential forms whose restriction to $H$ vanishes, and $A^*$ its annihilator,
\begin{eqnarray*}
A^* =\{\eta\,|\,\eta\wedge\theta=0 \quad\textrm{for all} \quad\theta\in\Theta^*\}.
\end{eqnarray*}
Elements of $A^*$ are called \emph{horizontal forms}.
\end{defi}
Choose locally (globally if $H$ admits a transverse orientation) a smooth $n-h$-form $\phi$ which is locally the wedge product of 1-forms from a basis of $\Theta^1$. Then $A^*$ consists of forms which are multiples of $\phi$.

\subsection{Existence of closed horizontal forms}

$k$-wealth, which subsumes existence of (at least local) submersions with horizontal fibers, implies abundance of closed horizontal $n-k$-forms, $n=dim(M)$. 

The abundance of closed horizontal $n-1$-forms can be seen in a direct manner too.

\begin{prop}
\label{n-1forms}
\emph{(\cite{Gromov-CC} page 156).}
On an $n$-dimensional Carnot manifold, every closed $n-1$-form is cohomologous to a (closed) horizontal form.
\end{prop}

\proof
The filtration of $H\subset H^{2}\subset \cdots\subset TM$ induces a filtration $A^{n-1}=\mathcal{F}^{1}\subset \mathcal{F}^{2}\subset \cdots\subset \Omega^{n-1}$ as follows : $\alpha\in\mathcal{F}^j$ if and only if there exist an $n$-form $\omega$ and a vectorfield $Z\subset H^j$ such that $\alpha=\iota_Z \omega$. We show that for all $j\geq 1$, $\mathcal{F}^{j+1} \subset \mathcal{F}^{j}+im(d)$. 

Let $\omega$ be an $n$-form and $X$, $Y$ be vectorfields such that $X\in H$, $Y\in H^{j}$. Then, using Lie derivatives,
\begin{eqnarray*}
\mathcal{L}_{X}(\iota_{Y}\omega)
&=&\iota_{\mathcal{L}_{X}(Y)} \omega+\iota_{Y}(\mathcal{L}_{X}\omega)\\
&=&\iota_{[X,Y]}\omega \textrm{ mod }\mathcal{F}^j .
\end{eqnarray*}
Thanks to Cartan's formula,
\begin{eqnarray*}
\mathcal{L}_{X}(\iota_{Y}\omega)
= d(\iota_{X}\iota_{Y}\omega)+\iota_{X}d(\iota_{Y}\omega)\in im(d)+\mathcal{F}^j .
\end{eqnarray*}
Therefore $\iota_{[X,Y]}\omega\in im(d)+\mathcal{F}^j$.

Let $\alpha\in \mathcal{F}^{j+1}$, $\alpha=\iota_Z \omega$ with $Z\in H^{j+1}$. Write $Z=\sum_\ell a_\ell [X_\ell ,Y_\ell ]$ where $a_\ell$ are functions, $X_\ell$, $Y_\ell$ are vectorfields, $X_\ell \in H$, $Y_\ell \in H^{j}$. Then $\alpha=\sum_\ell \iota_{[X_\ell ,Y_\ell ]}\omega_\ell$ (where $\omega_\ell =a_\ell \omega$), therefore $\alpha\in im(d)+\mathcal{F}^j$. This shows that $\mathcal{F}^{j+1} \subset \mathcal{F}^{j}+im(d)$. 

The bracket generating assumption, $H^r =TM$, implies that $\Omega^{n-1}=\mathcal{F}^{r}\subset im(d)+\mathcal{F}^{1}=im(d)+A^{n-1}$.
Given a closed $n-1$-form $\alpha$, the equation 
\begin{eqnarray*}
d\beta=-\alpha \quad\textrm{mod}\quad A^{n-1}
\end{eqnarray*}
admits a smooth solution $\beta\in \Omega^{n-2}$. Then $\alpha+d\beta\in A^{n-1}$ is a horizontal form.\qed

\subsection{A second proof of the isoperimetric inequality}

We prove once more that, for a Carnot group $G$ with Hausdorff dimension $Q$, bounded domains $D$ with piecewise smooth boundary satisfy
\begin{eqnarray*}
vol(D)\leq\textrm{const.}\,area(\partial D)^{Q/Q-1},
\end{eqnarray*}
where \emph{volume} (resp. \emph{area}) denotes $Q$-dimensional (resp. $Q-1$-dimensional) spherical Hausdorff measure.

\begin{rem}
\label{idim}
This implies Theorem \ref{thmi}. Indeed, the isoperimetric inequality for piecewise smooth domains always extends (possibly with a loss on constants) to arbitrary open sets.
\end{rem}

This second proof, which occupies the 4 next paragraphs, is borrowed from \cite{Gromov-CC} pages 167-168. It transposes a Euclidean argument which can be found for instance in \cite{Santalo}. It relies on homogeneity, scale invariance, an integration by parts and rearrangement. Horizontal differential forms, which have the right scale invariance under homothetic automorphisms, play a crucial role.

\subsection{Fundamental solution of the exterior differential} 

\begin{lemme}
\label{fundamental}
Let $G$ be a Carnot group of dimension $n$ and Hausdorff dimension $Q$, equipped with a left-invariant Carnot-Caratheodory metric. Let $p\in G$. There exists a smooth closed horizontal $n-1$-form $\omega_p$ on $G\setminus\{p\}$ such that
\begin{itemize}
  \item If $D$ is a bounded domain with piecewise smooth boundary and $p\in D$, 
\begin{eqnarray*}
\int_{\partial D}\omega_p =1. 
\end{eqnarray*}
  \item Write $\omega_p =\iota_X vol$ for some horizontal vectorfield $X$. Then
\begin{eqnarray*}
|X(q)|\leq\textrm{const.}\,|p-q|^{1-Q},
\end{eqnarray*}where $|p-q|$ is the Carnot-Caratheodory distance from $p$ to $q$.
\end{itemize}
\end{lemme}

\proof
Let $\delta_{\epsilon}$ denote the 1-parameter group of homothetic automorphisms of $(G,H=V^1)$ (see Definition \ref{defcg}). Consider the discrete group $Z=\{\delta_{2^n}\,|\,n\in\Z\}$. It acts properly discontinuously and cocompactly on $G\setminus\{e\}$, and preserves $H$. Therefore the quotient space $M=(G\setminus\{e\})/Z$ is a compact Carnot manifold without boundary. 

Let $S\subset G$ be a small (Euclidean) sphere centered at $e$. The map
\begin{eqnarray*}
S\times(0,+\infty)\to G\setminus\{e\},\quad
(q,\epsilon)\mapsto \delta_{\epsilon}(q)
\end{eqnarray*}
is a $Z$-equivariant diffeomorphism, where $Z$ acts trivially on $S$ and multiplicatively on $(0,+\infty)$. Therefore $M$ is diffeomorphic to $S \times S^{n-1}$. The cohomology class $c$, Poincar\'e-dual to the homology class of the $S^1$ factor, is represented by closed differential forms whose integral on the $S$ factor is equal to one. According to Proposition \ref{n-1forms}, one can choose a horizontal representative $\alpha$. Our $\omega_e$ is the pull-back of $\alpha$ under the covering map $\pi:G\setminus\{e\}\to M$, $\omega=\pi^* \alpha$, and $\omega_p$ is obtained from $\omega_e$ by left translation.

If $D$ is a domain that contains the ball $\beta$ bounded by $S$, then, by Stokes theorem,
\begin{eqnarray*}
\int_{\partial D}\omega_e -\int_{S}\omega_e =\int_{D\setminus \beta}d\omega_e =0.
\end{eqnarray*}
Thus 
\begin{eqnarray*}
\int_{\partial D}\omega_e 
&=&\int_{S}\pi^*\alpha \\
&=&\int_{\pi(S)}\alpha \\
&=& 1,
\end{eqnarray*}
by construction.

For $\epsilon=2^n$, $n\in\Z$, $\delta_{\epsilon}^{*}\omega_e =\omega_e$. From $\delta_{\epsilon}^{*}vol =\epsilon^Q vol$,  we get $(\delta_{\epsilon})_{*}X =\epsilon^{Q}X$, i.e. $X(\delta_{\epsilon}(q))=\epsilon^{Q-1}X(q)$, since $X$ is horizontal. In particular, the function
\begin{eqnarray*}
q\mapsto |q-e|^{Q-1}|X(q)|
\end{eqnarray*}
is invariant under the group $Z$, descends to a continuous function on the compact manifold $M$, so is bounded.\qed
 
\subsection{Integration versus area} 

\begin{lemme}
\label{area}
Let $\alpha=\iota_X vol$ be a horizontal $n-1$-form on $G$. Let $W\subset G$ be a hypersurface. Then
\begin{eqnarray*}
\int_{W}\alpha\leq\textrm{const.}\,\int_{W}|X(q)|\,dq
\end{eqnarray*}
where $dq$ denotes area, i.e. $Q-1$-dimensional spherical Hausdorff measure.
\end{lemme}

\proof
It suffices to verify the inequality for very small pieces of $W$, like intersections with small balls $B(q,\epsilon)$ centered on $W$. To save notation, let $q=e$. By definition of Hausdorff measure,
\begin{eqnarray*}
\epsilon^{1-Q}\int_{W\cap B(e,\epsilon)}dq \to 1
\end{eqnarray*}
as $\epsilon$ tends to 0. On the other hand,
\begin{eqnarray*}
\epsilon^{1-Q}\int_{W\cap B(e,\epsilon)}\alpha
&=&
\int_{\delta_{1/\epsilon}(W)\cap B(e,1)}\epsilon^{1-Q}\delta_{\epsilon}^{*}\alpha\\
&\to&\int_{V\cap \exp^{-1}B(e,1)}\iota_{X(e)}vol\\
&=&|X(e)|\,F(V,\frac{X(e)}{|X(e)|}),
\end{eqnarray*}
where $V$ is a hyperplane in the Lie algebra of $G$ that contains $V^2 \oplus\cdots\oplus V^r$. Since $F$ is a continuous function on a product of two projective spaces, it is bounded.\qed

\subsection{Integration by parts}

Let $D\subset G$ be a bounded domain with piecewise smooth boundary. Then
\begin{eqnarray*}
vol(D)&=&\int_{D}(\int_{\partial D}\omega_p )\,dp\\
&\leq&\int_{D\times\partial D}\textrm{const.}\,|p-q|^{1-Q}\,dq\,dp\\
&=&\textrm{const.}\,\int_{\partial D}(\int_{D}|p-q|^{1-Q}\,dp)\,dq.
\end{eqnarray*}

\subsection{Rearrangement} 

Let $B=B(q,R)$ be the Carnot ball centered at $q$ such that $vol(B)=vol(D)$. Then
\begin{eqnarray*}
\int_{D}|p-q|^{1-Q}\,dp
&=&\int_{D\cap B}|p-q|^{1-Q}\,dp+\int_{D\setminus B}|p-q|^{1-Q}\,dp\\
&\leq&\int_{D\cap B}|p-q|^{1-Q}\,dp+R^{1-Q}vol(D\setminus B)\\
&=&\int_{D\cap B}|p-q|^{1-Q}\,dp+R^{1-Q}vol(B\setminus D)\\
&\leq&\int_{D\cap B}|p-q|^{1-Q}\,dp+\int_{B\setminus D}|p-q|^{1-Q}\,dp\\
&=&\int_{B}|p-q|^{1-Q}\,dp\\
&=&\textrm{const.}\,R\\
&=&\textrm{const.}\,vol(D)^{1/Q}.
\end{eqnarray*}
Putting things together yields
\begin{eqnarray*}
vol(D)\leq\textrm{const.}\,vol(D)^{1/Q}\,area(\partial D),
\end{eqnarray*}
as expected.\qed

\section{The weight filtration of differential forms}

Gromov's integral geometric proof of the isoperimetric inequality in the previous section shows that horizontal $n-1$-forms can usefully replace sprays of horizontal curves. In this section, lower degree forms will be used as well in a H\"older exponent estimate as a replacement for foliations by higher dimensional horizontal submanifolds. Horizontal differential forms on a Carnot group are those which are contracted the most under homothetic automorphisms. This leads to the notion of weight. Gromov extracts a metric invariant from such weights. This yields upper bounds for the H\"older exponent $\alpha(M,H)$.

\subsection{Weights of differential forms}

Let $G$ be a Carnot group with Lie algebra $\mathcal{G}$. Left-invariant differential forms on $G$ split into homogeneous components under the homothetic automorphisms $\delta_\epsilon$,
\begin{eqnarray*}
\Lambda^* \mathcal{G}^* =\bigoplus_{w}\Lambda^{*,w} \quad \textrm{where}\quad \Lambda^{*,w} =\{\alpha\,|\,\delta_{\epsilon}^{*}\alpha=\epsilon^{w}\alpha\}.
\end{eqnarray*}

\begin{exemple}
\label{weightheis}
If $G=Heis^{2m+1}$ is the Heisenberg group, for each degree $q\not=0$, $2m+1$,
\begin{eqnarray*}
\Lambda^q \mathcal{G}^* =\Lambda^{q,q}\oplus\Lambda^{q,q+1},
\end{eqnarray*}
where $\Lambda^{q,q}=\Lambda^{q}(V^1 )^*$ and $\Lambda^{q,q+1}=\Lambda^{q-1}(V^1 )^* \otimes (V^2 )^*$.
\end{exemple}

This gradation by weight depends on the group structure. What remains for general Carnot manifolds is a filtration.

\begin{defi}
\label{defweight}
Let $(M,H)$ be a Carnot manifold, $m\in M$. Say a $q$-form $\alpha$ on $T_m M$ has \emph{weight} $\geq w$ if it vanishes on $q$-vectors of $H^{i_1}\otimes\cdots \otimes H^{i_q}$ whenever $i_1 + \cdots +i_q <w$. If $(M,H)$ is equiregular, such forms constitute a subbundle $\Lambda^{q,\geq w}T^* M$. The space of its smooth sections is \emph{denoted by} $\Omega^{*,\geq w}$.
\end{defi}
Note that each $\Omega^{*,\geq w}$ is a differential ideal in $\Omega^*$. 

\begin{exemple}
\label{exweighthor}
Assume $(M,H)$ is equiregular of dimension $n$ and Hausdorff dimension $Q$. Then a differential $q$-form on $M$ is horizontal if and only if it has maximal weight, i.e. weight $\geq Q-n+q$.
\end{exemple}

\begin{lemme}
\label{grad}
Assume $(M,H)$ is equiregular. The graded algebra $\bigoplus_{w}\Lambda^{*,\geq w}/\Lambda^{*,\geq w+1}$ identifies with the space $\Lambda^* \mathcal{G}_m$ of left invariant differential forms on the tangent Carnot group $G_m$. 
\end{lemme}

\proof
According to \cite{NSW}, the Lie algebra $\mathcal{G}_m$ is the graded space 
\begin{eqnarray*}
\mathcal{G}_m = gr(H_{m}^{\bullet})=\bigoplus_{i=1}^{r} H_{m}^{i}/H_{m}^{i-1}
\end{eqnarray*}
associated to the filtration $(H_{m}^{i})_{1\leq i\leq r}$ of the tangent space $T_m M$, equipped with a bracket induced by the Lie bracket on vectorfields. Since the filtration of differential forms is defined by duality,
\begin{eqnarray*}
gr(\Lambda^{*,\bullet})=gr(\Lambda^* (H_{m}^{\bullet})^*) )=\Lambda^* (gr(H^{\bullet}))^* =\Lambda^* (\mathcal{G}_m)^* .\qed
\end{eqnarray*}

\begin{cor}
\label{integerweights}
The values taken by weights for an equiregular Carnot manifold are those of Carnot groups, i.e. integers between 1 and $Q$.
\end{cor}

\begin{defi}
\label{defweightinv}
Let $(M,H)$ be a Carnot manifold. Define the \emph{weight invariant} $W_q (M,H)$ as the largest $w$ such that there exists arbitrarily small open sets with smooth boundary $U\subset M$ and nonzero classes in $H^q (U,\R)$ which can be represented by closed differential forms of weight $\geq w$.
\end{defi}

\begin{exemple}
\label{exweightinv}
By Proposition \ref{n-1forms}, all Carnot manifolds of dimension $n$ and Hausdorff dimension $Q$ have $W_{n-1}\geq Q-1$.
\end{exemple}

\subsection{Straight cochains}

Following Gromov (\cite{Gromov-CC}, pages 247-249), we define a metric invariant whose behaviour is similar to $W_q$. The starting point is the following characterization of weight on Carnot groups.

\begin{rem}
\label{remweightcontr}
A differential form $\omega$ on a Carnot group has weight $\geq w$ if and only if 
\begin{eqnarray*}
\n{\delta_{\epsilon}^{*}\omega}_{\infty}\leq \textrm{const.}\,\epsilon^w \quad\textrm{for}\quad \epsilon\leq 1.
\end{eqnarray*}
\end{rem}

The next idea is to replace differential forms, as a tool for cohomology calculations, by Alexander-Spanier straight cochains, which have the advantage of being functorial under homeomorphisms.

\begin{defi}
\label{defas}
\emph{(Alexander-Spanier).} Let $X$ be a metric space, and $t>0$. A \emph{straight $q$-cochain of size $t$} on $X$ is a bounded function on $q+1$-tuples of points of $X$ of diameter less than $t$. The \emph{$\epsilon$-absolute value} a straight $q$-cochain $c$ is its $\ell^{\infty}$ norm as a cochain of size $\epsilon$, i.e.
\begin{eqnarray*}
|c|_{\epsilon}=\sup\{|c(\sigma)|\,|\,q+1-\textrm{tuples }\sigma\textrm{ of diameter }<\epsilon\}.
\end{eqnarray*}
\end{defi}
Straight cochains of size $t$ form a complex, since they coincide with simplicial cochains on the simplicial complex $X_t$ with vertex set $X$, such that $q+1$ vertices span a $q$-simplex if and only if their mutual distances in $X$ are less than $t$. The simplicial chains on $X_t$ are called straight chains of size $t$.

If $X$ is a compact manifold with boundary, or biH\"older homeomorphic to such, then, for $t$ small enough, straight chains (resp. cochains) of size $t$ compute homology (resp. cohomology). Given a cohomology class $\kappa$ and a number $\nu>0$, one can define the \emph{$\nu$-norm}
\begin{eqnarray*}
\n{\kappa}_{\nu}=\liminf_{\epsilon\to 0}\epsilon^{-\nu}\inf\{|c|_{\epsilon}\,|\,\textrm{ cochains } c \textrm{ of size } \epsilon \textrm{ representing } \kappa\}.
\end{eqnarray*}

The next two propositions provide opposite estimates on norms.

\begin{prop}
\label{lowernorm}
In a Riemannian manifold with boundary, all straight cocycles $c$ representing a nonzero class $\kappa$ of degree $q$ satisfy
\begin{eqnarray*}
|c|_{\epsilon}\geq \textrm{const.}(\kappa)\,\epsilon^q .
\end{eqnarray*}
In other words, $\n{\kappa}_{q}>0$.
\end{prop}

\proof
Fix a cycle $c'$ such that $\kappa(c')>0$. Subdivide it as follows : fill simplices with geodesic singular simplices, subdivide them and keep only their vertices. This does not change the homology class. The number of simplices of size $\epsilon$ thus generated is $\leq \textrm{const.}(c')\,\epsilon^{-q}$. For any representative $c$ of size $\epsilon$ of $\kappa$,
\begin{eqnarray*}
\kappa(c')=c(c')\leq \textrm{const.}\,\epsilon^{-q} |c|_{\epsilon}.\qed
\end{eqnarray*}

\begin{prop}
\label{uppernorm}
Let $(M,H)$ be an equiregular Carnot manifold. Let $U\subset M$ be a bounded open set with smooth boundary. Let $\omega$ be a closed differential form on $U'$ of weight $\geq w$. Then, for every $\epsilon$ small enough, the cohomology class $\kappa\in H^q (U,\R)$ of $\omega$ can be represented by a straight cocycle $c_{\epsilon}$ (maybe defined on a slightly smaller homotopy equivalent open set) such that 
\begin{eqnarray*}
|c_{\epsilon}|_{\epsilon}\leq \textrm{const.}\,\epsilon^w .
\end{eqnarray*}
In other words, $\n{\kappa}_{w}<+\infty$.
\end{prop}

\proof
In the case of a Carnot group $G$. Use simultaneously left-invariant Riemannian and Carnot-Caratheodory metrics. Use the exponential map to push affine simplices in the Lie algebra to the group. Fill in all straight simplices in $G$ of unit Carnot-Caratheodory size with such affine singular simplices. Note that the Carnot-Caratheodory diameters and the Riemannian volumes of these singular simplices ($C^1$ maps of the standard simplex to $G$) are bounded by some constant $V$. Then forget th Riemannian metric, apply $\delta_{\epsilon}$ and obtain a filling $\sigma_{\epsilon}$ for each straight simplex $\sigma$ in $G$ of Carnot-Caratheodory size $\epsilon$. Let $U'$ be an open set whose closure is contained in $U$ and which is a deformation retract of $U$. For $\epsilon$ small enough, the filling of a straight simplex of size $\epsilon$ in $U'$ is contained in $U$. Define a straight cochain $c_{\epsilon}$ of size $\epsilon$ on $U'$ by
\begin{eqnarray*}
c_{\epsilon}(\sigma)=\int_{\sigma_{\epsilon}}\omega.
\end{eqnarray*}
Since $\omega$ is closed, Stokes theorem shows that $c_{\epsilon}$ is a cocycle. Its cohomology class in $H^q (U',\R)\simeq H^q (U,\R)$ is the same as $\omega$'s. Furthermore,
\begin{eqnarray*}
|c_{\epsilon}(\sigma)|
&=&\int_{\sigma_1}\delta_{\epsilon}^{*}\omega\\
&\leq&V\,\n{\delta_{\epsilon}^{*}\omega}_{\infty}\\
&\leq&\textrm{const.}(\omega)\,\epsilon^{w}.\qed
\end{eqnarray*}

\subsection{The metric weight invariant}

Here is the promised metric analogue for $W_q$.

\begin{defi}
\label{defmetweightinv}
Let $X$ be a metric space. Define $MW_q (X)$ as the supremum of numbers $\nu$ such that there exist arbitrarily small open sets $U\subset M$ and nonzero straight cohomology classes $\kappa\in H^q (U,\R)$ with finite $\nu$-norm $\n{\kappa}_{\nu}<+\infty$.
\end{defi}

\begin{exemple}
\label{exmw}
Proposition \ref{lowernorm} shows that Euclidean space has $MW_q \leq q$.

Proposition \ref{uppernorm} shows that equiregular Carnot manifolds satisfy $MW_q \geq W_q$. 
\end{exemple}

\begin{prop}
\label{holderweight}
Let $f:X\to Y$ be a $C^{\alpha}$-H\"older continuous homeomorphism. Let $\kappa\in H^q (Y,\R)$. Then
\begin{eqnarray*}
\n{\kappa}_{\nu}<+\infty\Rightarrow\n{f^{*}\kappa}_{\nu\alpha}<+\infty.
\end{eqnarray*}
In particular, $MW_q (X)\geq \alpha MW_q (Y)$.
\end{prop}

\proof
If $\sigma$ is a straight simplex of size $\epsilon$ in $X$, $f(\sigma)$ has size $\epsilon' \leq \n{f}_{C^{\alpha}}\,\epsilon^{\alpha}$ in $Y$. If $c$ is a representative of $\kappa$, $f^* c$ is a representative of  $f^{*}\kappa$, and
\begin{eqnarray*}
\epsilon'^{-\nu}|c|_{\epsilon'}
&\geq&\epsilon'^{-\nu}|c(f(\sigma))|\\
&=&\epsilon'^{-\nu}|f^* c(\sigma)|\\
&\geq&\n{f}_{C^{\alpha}}^{-\nu}\,\epsilon^{-\nu\alpha}|f^* c(\sigma)|.
\end{eqnarray*}
Therefore
\begin{eqnarray*}
\epsilon^{-\nu\alpha}|f^* c|_{\epsilon}
\leq \n{f}_{C^{\alpha}}^{\nu}\,\epsilon'^{-\nu}|c|_{\epsilon'}.
\end{eqnarray*}
This leads to
\begin{eqnarray*}
\n{f^* \kappa}_{\nu\alpha}\leq \n{f}_{C^{\alpha}}^{\nu}\,\n{\kappa}_{\nu}.\qed
\end{eqnarray*}

\begin{cor}
\label{corweighthol}
Let $(M,H)$ be an equiregular Carnot manifold. Then for all $q=1,\ldots,n-1$, 
\begin{eqnarray*}
\alpha(M,H)\leq \frac{q}{W_q (M,H)}.
\end{eqnarray*}
\end{cor}

\section{Complexes of differential forms}

We want that every closed form be cohomologous to another one of high weight. In other words, we need to compute cohomology with a subcomplex of differential forms of rather high weights. Such complexes have been discovered by M. Rumin, first for Heisenberg groups, \cite{Rumin-94}, later on for general Carnot groups, \cite{Rumin-99}, and equihomological equiregular Carnot manifolds, \cite{Rumin-05}.

\subsection{Rumin's contact complex}

We begin with the special case of contact manifolds. 

 \begin{defi}
\label{defI}
Let $(M,H)$ be a contact manifold, i.e. $H=ker(\theta)$ where $\theta$ is an ordinary 1-form, and $d\theta_{|H}$ is symplectic. Let
\begin{eqnarray*}
I^* =\{\alpha\wedge\theta+\beta\wedge d\theta\,|\,\alpha,\, \beta\in\Omega^*\}
\end{eqnarray*}
denote the differential ideal generated by $\theta$. Let
\begin{eqnarray*}
J^* =\{\eta\in\Omega^* \,|\,\eta\wedge\gamma=0 \quad\textrm{for all}\quad\gamma\in I^* \}
\end{eqnarray*}
denote its annihilator.
\end{defi}
Let us denote by $d^H :I^* \to I^*$ (resp. $d_H :J^* \to J^*$) the operators induced by the exterior differential. 

\begin{theo}
\label{rumin94}
\emph{(M. Rumin, \cite{Rumin-94}).} Let $(M,H)$ be a $2m+1$-dimensional contact manifold. There exists a second order differential operator $D:\Omega^m /I^m \to J^{m+1}$ such that the complex
\begin{eqnarray*}
0\to \Omega^1 /I^1 \fleche{\dh}\cdots\fleche{\dh} \Omega^m /I^m \fleche{D} J^{m+1} \fleche{d_{H}}\cdots\fleche{d_{H}} J^{2m+1}\to 0
\end{eqnarray*}
is homotopy equivalent to the de Rham complex. In particular, it computes the cohomology of $M$.
\end{theo}

\proof
Start with the 3-dimensional case. View elements $\eta$ of the quotient space $\Omega^1 /I^1$ as partially defined 1-form (i.e. sections of the dual bundle $H^*$). When is such a partial form the restriction of the differential of a function ? There should exist an extension $\alpha$ of $\eta$ in the missing direction $Z$ which is a closed 1-form. Write $Z=[X,Y]$ where $X$ and $Y$ are horizontal vectorfields. Since we want in particular that $d\alpha(X,Y)=0$, there is only one choice for $\alpha(Z)$,
\begin{eqnarray*}
\alpha(Z)=\alpha([X,Y])=X\eta(Y)-Y\eta(X).
\end{eqnarray*} 
Thus the necessary condition for $\eta$ to be a $d^{H}u$ is $d\alpha=0$. Since $\alpha$ depends on first derivatives of $\eta$, $D\eta=d\alpha$ depends on second derivatives. By construction, $D\eta(X,Y)=0$, i.e. the 2-form $D\eta$ is horizontal, $D\eta\in J^{2}$. 

Locally, a closed form is exact, so $D\eta=0$ implies $\eta=d^{H}u$ locally. If $\gamma$ is a closed horizontal 2-form, then locally $\gamma=d\alpha=D(\alpha_{|H})$. In other words, the sequence of differential operators 
\begin{eqnarray*}
\Omega^0 \fleche{d^{H}} \Omega^1 /I^1 \fleche{D} J^{2} \fleche{d_{H}} J^{3},
\end{eqnarray*}
called the \emph{Rumin complex}, is locally exact. Therefore, globally, it is as good as the full de Rham complex for cohomology calculations. A homotopy equivalence is a sharper way to express this. A \emph{homotopy equivalence of complexes} $C^*$ and $D^*$ is a pair of maps $h:C^* \to D^*$ and $h':D^* \to C^*$ satisfying of course $dh=hd$, $dh'=h'd$, such that there exist $B:C^* \to C^*$ and $B':D^* \to D^*$ such that $h'h=1-dB-Bd$, $hh'=1-dB'-B'd$. The Rumin complex trivially injects into the de Rham complex, except in degree 1. In that degree, we set $h(\eta)=\alpha$, the unique extension just described. Conversely, the de Rham complex trivially maps to our complex except in degree 2. If $\alpha$ is a 1-form, $(hh'-1)(\alpha)(Z)= X\eta(Y)-Y\eta(X)-\alpha([X,Y])=d\alpha(X,Y)$. Thus, assuming that $\theta(Z)=1$, we define $B':\Omega^{2}\to\Omega^{1}$ by $B'(\gamma)=-\gamma(X,Y)\theta$, and decide that $B'=0$ in other degrees. If $\gamma$ is a 2-form, the wished identity $(hh'-1)(\gamma)=-dB'(\gamma)$ suggests to define $h'(\gamma)=\gamma+d(\gamma(X,Y)\theta)$, which indeed belongs to $J^{2}$. With this choice, $hh'=1-dB'-B'd$. On the Rumin complex, one simply takes $B=0$. Since $h'h=1$, everything fits nicely.

The higher dimensional case requires more care, and we only explain the construction of $D$, referring to Rumin's original paper for the homotopy equivalence. Let $\eta\in\Omega^m$. What is the condition for $\eta$ mod $I^m$ to be in the image of $d^H :\Omega^{m-1}/I^{m-1}\to \Omega^m /I^m$ ? There should exist $\gamma=\alpha\wedge\theta+\beta\wedge d\theta\in I^m$ such that $d(\eta+\gamma)=0$. Note that
\begin{eqnarray*}
d\gamma=d(\alpha\wedge\theta)+d\beta\wedge d\theta=d((\alpha+(-1)^m d\beta)\wedge\theta),
\end{eqnarray*}
so that $\beta$ does not bring anything new. Compute
\begin{eqnarray*}
d(\eta+\alpha\wedge\theta)=d\eta+d\alpha\wedge\theta+(-1)^{m-1}\alpha\wedge d\theta.
\end{eqnarray*}
Since $d\theta_{|H}$ is symplectic, wedging with $d\theta_{|H}:\Lambda^{m-1}H^* \to\Lambda^{m+1}H^*$ is a bijection. Therefore there is a unique $\alpha_{|H}$ (and therefore a unique $\alpha\wedge\theta$) such that
\begin{eqnarray*}
\alpha_{|H}\wedge d\theta_{|H}=-d\eta_{|H}.
\end{eqnarray*}

Let us denote $d(\eta+\alpha\wedge\theta)$ by $D\eta$. By construction, $D\eta_{|H}=0$, thus $D\eta\wedge\theta=0$. Furthermore, there exists a form $\kappa$ such that $D\eta=\kappa\wedge\theta$. Then $0=d(D\eta)=d\kappa\wedge\theta+(-1)^m \kappa\wedge d\theta$, showing that $\kappa_{|H}\wedge d\theta_{|H}=0$. This implies that $D\eta\wedge d\theta=0$, i.e. $D\eta\in J^{m+1}$. Thus $D$ connects the complexes $d^H :I^* \to I^*$ and $d_H :J^* \to J^*$ into a single complex.\qed

\begin{cor}
\label{weightcontact}
If $(M,H)$ is a $2m+1$-dimensional contact manifold, then $W_q (M,H)\geq q+1$ for all $q\geq m+1$.
\end{cor}

\begin{rem}
\label{reminvcontact}
Rumin's contact complex does not involve any arbitrary choices, it is invariant under all contactomorphisms. However, the maps $h'$ and $B'$ depend on the choice of a contact form $\theta$ or a complement $Z$, $Z\oplus H=TM$.
\end{rem}

\subsection{The weight preserving part of $d$}

The fact that $B=0$ in Rumin's contact complex suggests to view it as the subcomplex $im(h)$ in the de Rham complex, and $h'$ as a retraction. In spite of its cost (loss of invariance), this point of view opens the way to generalization.

We are looking for a retraction $r=1-dB'-B'd$ of the de Rham complex, onto forms of high weight. In other words, we want $r$ to kill low weight components of forms. But $B'$ inverts the exterior differential, inasmuch as possible. Thus $B'$ should especially invert $d$ on low weights. How does $d$ behave with weights ?

\begin{lemme}
\label{dzero}
Each $\Omega^{*,\geq w}$ is a differential ideal in $\Omega^*$. The operator $d^0$ induced on \begin{eqnarray*}
C^{\infty}(\bigoplus_{w}\Lambda^{*,\geq w}/\Lambda^{*,\geq w+1})=\bigoplus_{w}\Omega^{*,\geq w}/\Omega^{*,\geq w+1}
\end{eqnarray*}
from the exterior differential on $\Omega^{*}$ is algebraic (it does not depend on derivatives), and acts fiberwise as the exterior differential on left-invariant forms.
\end{lemme}

\proof
If $\omega$ is a left invariant form of weight $\geq w$ and $f$ a function, $df\wedge\omega$ has weight $\geq w+1$, so $d^0 (f\omega)=fd^0 \omega$, and $d^0$ is algebraic. The pointwise computation of $d^0$ is done in the Carnot group case. There, since $d$ commutes with pull-back under homothetic automorphisms, if $\phi$ is a left-invariant homogeneous differential form of weight $w$, so is $d\phi$. This shows that $d^0 =d$ on left invariant forms.\qed

\begin{defi}
\label{gradedhomology}
On the tangent Carnot Lie algebra $\mathcal{G}_m$, $d^0$ commutes with the homothetic automorphisms $\delta_{\epsilon}$, thus its cohomology is graded by degree and weight,
\begin{eqnarray*}
H^* (\mathcal{G}_m ,\R)=\bigoplus_{q,w}H^{q,w}(\mathcal{G}_m ,\R).
\end{eqnarray*}
\end{defi}

\begin{exemple}
\label{dzeroheis}
Cohomology of the Heisenberg Lie algebra $\mathcal{H}^{2m+1}$.
\end{exemple}
Choose $\theta$ in $(V^{2})^*$. Every $\alpha\in\Lambda^* \mathcal{G}^*$ can be uniquely written $\alpha=\eta +\theta\wedge\beta$ with $\eta$, $\beta\in \Lambda^* (V^1 )^*$. Then
\begin{eqnarray*}
d^0 \alpha=d\theta\wedge\beta=:L\beta,
\end{eqnarray*}
where $d\theta$ is a symplectic form on $V^1$. Thus cohomology splits into  $H^{q,q+1}=\theta \wedge ker(L)$ and $H^{q,q}=\Lambda^q (V^1 )^* /im(L)$.

\subsection{Rumin's complex for equihomological Carnot manifolds}

Lemma \ref{dzero} suggests that $B'$ should be an inverse of $d^0$. This operator is defined on a quotient $\bigoplus_{w}\Omega^{*,\geq w}/\Omega^{*,\geq w+1}$. If $(M,H)$ is equiregular, choices of complements $V^k$ of $H^{k-1}$ in $H^k$ allow a lift of $d^0$ to an operator on differential forms.

\begin{exemple}
\label{dzerocontact}
Let $(M,H)$ be a contact manifold. A choice of contact form $\theta$ determines a complement $V^2 =ker(d\theta)$. Every differential form can be uniquely written $\alpha=\eta +\theta\wedge\beta$ with $\iota_{V^2}\eta=\iota_{V^2}\beta=0$. Then
\begin{eqnarray*}
d^0 \alpha=d\theta\wedge\beta.
\end{eqnarray*}
\end{exemple}

\begin{defi}
\label{defequihom}
Say an equiregular Carnot manifold is \emph{equihomological} if the dimensions of the cohomology spaces of tangent Lie algebras and of their weight filtrations are constant.
\end{defi}

When this is the case (for instance, for Carnot groups), one can smoothly choose complements
\begin{itemize}
  \item $V_j$ of $H^{j-1}$ in $H^j$,
  \item $F$ of $ker(d^0 )$ in $\Lambda^* T^* M$,
  \item $E$ of $im(d^0 )$ in $\ker(d_0 )$.
\end{itemize} 
The choice of $V_j$ allows to view $d^0$ as acting on forms (and not on some quotient space) and so the other choices make sense in turn.

This determines an inverse $(d^{0})^{-1}$, with kernel $E+F$ and image $F$. Set
\begin{eqnarray*}
r=1-d(d^{0})^{-1}-(d^{0})^{-1}d.
\end{eqnarray*}
$r$ is a first order differential operator, compatible with weight filtrations.

\begin{theo}
\label{rumgen}
\emph{(M. Rumin, \cite{Rumin-99}).}
Let $(M,H)$ be an equihomological equiregular Carnot manifold. The iterates $r^j$ stabilize to a projector $p$ of $\Omega^*$, with image the subcomplex
\begin{eqnarray*}
\mathcal{E}=ker((d^{0})^{-1})\cap ker((d^{0})^{-1}d)
=\{\eta\in C^{\infty}(E\oplus F)\,|\,d\eta\in C^{\infty}(E\oplus F)\}.
\end{eqnarray*}
Furthermore, if $\pi$ denotes the projector onto $E$ with kernel $im(d^0 )\oplus F$, then, on $\mathcal{E}$, $p\circ\pi=identity$. 
\end{theo}

\proof
See \cite{Rumin-99}.\qed

\begin{cor}
\label{corweight}
Let $(M,H)$ be an equihomological equiregular Carnot manifold. Assume that there exists a point $m\in M$ such that, in the cohomology of the tangent Lie algebra $\mathcal{G}_m$, $H^{q,w'}(\mathcal{G}_m ,\R)=0$ for all $w'<w$. Then $W_q (M,H)\geq w$, and, as a consequence, $\alpha(M,H)\leq q/w$.
\end{cor}

\proof
The vanishing of $H^{q,w'}(\mathcal{G},\R)$ is an open condition on a Lie algebra. Therefore the assumptions are satisfied in a neighborhood of $m$. On this neighborhood,
by equihomologicality, $dim(\bigoplus_{w'\geq w}H^{q,w'}(\mathcal{G}_{m'} ,\R))$ is constant, thus one can choose a smooth complement $E^q$ of $im(d^0)\cap\Lambda^{q,\geq w}$ in $ker(d^0 )\cap\Lambda^{q,\geq w}$, complete it into a complement $E$ of $im(d^0)$ in $ker(d^0 )$ and pick a smooth complement $F$ of $ker(d^0 )$ in $\Lambda^*TM$. Let $U$ be some smaller neighborhood of $m$ with smooth boundary, such that $H^q (U,\R)\not=0$. Let $\kappa$ be a nonzero class in $H^q (U,\R)$. By Theorem \ref{rumgen}, $\kappa$ contains a closed form $\omega$ which belongs pointwise to $E\oplus F$. Thus $\pi\omega\in E^q$ has weight $\geq w$. Since $p$ is weight-preserving, $\omega=p\circ\pi(\omega)$ has weight $\geq w$ too. This shows that $W_q (M,H)\geq w$. The conclusion $\alpha(M,H)\leq q/w$ then follows from Corollary \ref{corweighthol}.\qed

\subsection{Duality}

The weight gradation of Lie algebra cohomology is invariant under Poincar\'e duality. This is useful for calculating examples.

\begin{prop}
\label{duality}
Let $G$ be a Carnot group with dimension $n$ and Hausdorff dimension $Q$. Then, in the cohomology of its Lie algebra $\mathcal{G}$, $H^{q,w}(\mathcal{G},\R)$ is isomorphic to $H^{n-q,Q-w}(\mathcal{G},\R)$.
\end{prop}

\proof
Choose a Euclidean structure on $\mathcal{G}$ which makes all $V^j$ orthogonal. Observe that the corresponding Hodge $*$-operator maps $\Lambda^{q,w}$ to $\Lambda^{n-q,Q-w}$. The adjoint $\delta^0$ of $d^0$ is given by $\delta^0 =\pm *d^0 *$. Choose 
\begin{eqnarray*}
E=\textrm{ orthogonal complement of } im(d^0 ) \textrm{ in }ker(d^0 )=ker(\delta^0 )\cap ker(d^0 ).
\end{eqnarray*}
Then $E$ is graded, and $E^{q,w}$ maps isomorphicly to $H^{q,w}$ in cohomology. Since $*E^{q,w}=E^{n-q,Q-w}$, the conclusion follows.\qed

\subsection{Examples}

\textbf{Degree $n-1$}. On any Carnot group, the space of closed invariant 1-forms is $(V^1 )^* =\Lambda^{1,1}$, thus $H^{1,w}(\mathcal{G})=0$ for $w>1$. Proposition \ref{duality} implies that $H^{n-1,w}(\mathcal{G})=0$ for $w<Q-1$, and Corollary \ref{corweight} gives $W_{n-1}(G)\geq Q-1$. We already knew this from Proposition \ref{n-1forms}.

\medskip

\textbf{Contact case}. On a symplectic $2m$-space, $L$ is injective in degrees $\leq m-1$ and surjective onto degrees $\geq m+1$. Therefore, following Example \ref{dzeroheis}, for $2m+1$-dimensional contact manifolds, $H^{q,q}=0$ for $q\geq m+1$. This implies $W_q (M,H)\geq q+1$ for $q\geq m+1$, a fact we already knew from Corollary \ref{weightcontact}.

\medskip

\textbf{Quaternionic contact case}. The quaternionic Heisenberg Lie algebra is $\mathcal{G}=V^1 \oplus V^2$ where $V^1 =\mathbf{H}^m$, $V^2 =\Im m(\mathbf{H}) $ and for $X$, $Y\in V^1$, $[X,Y]=\Im m\langle X,Y\rangle$. The group $Sp(n)Sp(1)$ acts by automorphisms on $\mathcal{G}$, and
\begin{eqnarray*}
\Lambda^{2,*}=\Lambda^{2,2}\oplus\Lambda^{2,3}\oplus\Lambda^{2,4}
\end{eqnarray*}
is a decomposition into irreducible summands. $d^0$ does not vanish identicly on $\Lambda^{2,3}$ or on $\Lambda^{2,4}$. Therefore $d^0$ is injective on these subspaces, and $H^{2,3}(\mathcal{G},\R)=H^{2,4}(\mathcal{G},\R)=0$. This implies $W_{n-2}(G)\geq Q-2=4m+4$.

\medskip

\textbf{Rank 2 distributions}. If $\mathcal{G}$ is a Carnot Lie algebra with $dim(V^1 )=2$, then $V^2 =[V^1 ,V^1 ]$ is 1-dimensional and $[,]:\Lambda^{2}(V^1 ) \to V^2$ is injective. Its adjoint $d^0 :\Lambda^{1,2}=(V^2 )^*=\to \Lambda^{2}(V^1 )^* =\Lambda^{2,2}$ is onto, and $H^{2,2}=0$. Furthermore, if $dim(V^3 )\not=1$, $[,]:V^1 \otimes V^2 \to V^3$ is injective. Its adjoint $d^0 :\Lambda^{1,3}=(V^3 )^* \to (V^1 )^* \otimes (V^2 )^* =\Lambda^{2,3}$ is onto, and $H^{2,3}=0$. Hence, for any equihomological equiregular Carnot manifold $(M,H)$ with $dim(H)=2$, $W_{2}(M,H)\geq 3$, $\alpha(M,H)\leq 2/3$. If furthermore $dim(H^3 )\geq 5$, then $W_{2}(M,H)\geq 4$, $\alpha(M,H)\leq 1/2$. Note that this bound is always worse that what is obtained when considering $n-1$-forms.

\subsection{Back to regular isotropic subspaces}

\begin{prop}
\label{weightregiso}
Let $(M,H)$ be an equiregular Carnot manifold. If $H$ contains a regular isotropic horizontal $k$-plane at some point $m$, then $H^{k,w}(\mathcal{G}_m ,\R)=0$ for all $w\geq k+1$. If $(M,H)$ is furthermore equihomological, it follows that $W_{n-k}(M,H)\geq Q-k$ and $\alpha(M,H)\leq \frac{n-k}{Q-k}$.
\end{prop}

\proof
If $S\subset H_m$ is regular isotropic, then $S$ viewed as a subspace in $V^1 \subset \mathcal{G}_m$ is regular isotropic for $G_m$ as well. Thus is the sequel, $M=G_m$ is a Carnot group, $H=V^1$ is defined by a left invariant $\R^{n-h}$-valued 1-form $\theta=(\theta_1 ,\ldots,\theta_{n-h})$.

Let $w\geq k+1$ and $\omega\in ker(d^0 )\cap \Lambda^{k,w}$. Then there exists an $\R^{n-h}$-valued $k-1$-form $\eta=(\eta_1 ,\ldots,\eta_{n-h})$ such that
\begin{eqnarray*}
\omega=\sum_{i=1}^{n-h}\theta_i \wedge \eta_i =\theta\wedge\eta.
\end{eqnarray*}
Then 
\begin{eqnarray*}
(d^0 \omega)_{|V^1}=(d^0 \theta)_{|V^1}\wedge \eta_{|V^1}.
\end{eqnarray*}
Since $S$ is isotropic, for $X\in V^1$,
\begin{eqnarray*}
(\iota_X d^0 \omega)_{|S}=(\iota_X d^0 \theta)_{|S}\wedge \eta_{|S}.
\end{eqnarray*}
Choose a Euclidean structure on $S$ and use its Hodge $*$ operator. By regularity, one can choose a vector $X \in V^1$ such that $(\iota_{X}d^0 \theta)_{|S}=*(\eta_{|S})$. Then
\begin{eqnarray*}
(\iota_X d^0 \omega)_{|S}=*(\eta_{|S})\wedge \eta_{|S}
\end{eqnarray*}
is nonnegative. Since $d^0 \omega=0$, this implies $\eta_{|S}=0$. Since the variety of isotropic $k-1$-subspaces is smooth, the linear span of the set of decomposable $k-1$-vectors associated to isotropic subspaces is all of $\Lambda^{k-1} V^1$, so $\eta_{|V^1}=0$ and $\omega=0$. One concludes that $ker(d^0 )\cap \Lambda^{k,w}=0$. In particular, $H^{k,w}(\mathcal{G}_m ,\R)=0$. 

By Poincar\'e duality (Proposition \ref{duality}), $H^{n-k,w}(\mathcal{G}_m ,\R)=0$ for all $w<Q-k$. The conclusion follows from Corollary \ref{corweight}.\qed

\begin{rem}
\label{genericweight}
According to Proposition \ref{isoreggeneric}, Proposition \ref{weightregiso} applies to generic $h$-dimen\-sional distributions on $n$-dimensional manifolds, provided $h-k\geq (n-h)k$.
\end{rem}

\section{Conclusion}

Concerning the H\"older equivalence problem, the direct approach using differential forms seems to cover all results obtained by the horizontal submanifold method, and has a wider scope (see \cite{Rumin-05}). Nevertheless, the bounds obtained are never sharp, even in the case of the 3-dimensional Heisenberg group. New ideas are needed.

\par\medskip\noindent
Pierre Pansu\\
Laboratoire de Math\'ematique d'Orsay\\
UMR 8628 du C.N.R.S.\\
B\^atiment 425\\
Universit\'e Paris-Sud - 91405 Orsay (France)\\
\smallskip\noindent
{\tt\small Pierre.Pansu@math.u-psud.fr}\\
http://www.math.u-psud.fr/$\sim$pansu

\end{document}